\documentclass[aap,MSNbibl,citesort,dvips]{arximspdf}
\usepackage{graphicx}

\doi{10.1214/09-AAP673}
\volume{20}
\issue{5}
\pubyear{2010}
\firstpage{1801}
\lastpage{1830}

\makeatletter

\newtheorem{theorem}{Theorem}
\newtheorem{lemma}[theorem]{Lemma}
\newtheorem{proposition}[theorem]{Proposition}

\newproclaim{definition}{Definition}
\newproclaim{remark}[definition]{Remark}

\def\d{{d}}
\def\i{{i}}

\makeatother

\begin{document}
\begin{frontmatter}

\title{Wiener--Hopf factorization and distribution of extrema for a
family of L\'evy processes}
\runtitle{Wiener--Hopf factorization and distribution of extrema}

\begin{aug}
\author[A]{\fnms{Alexey} \snm{Kuznetsov}\corref{}\thanksref
{t1}\ead[label=e1]{kuznetsov@mathstat.yorku.ca}}
\runauthor{A. Kuznetsov}
\affiliation{York University}
\address[A]{Department of Mathematics and Statistics\\
York University \\
Toronto, Ontario, M3J 1P3\\
Canada\\
\printead{e1}} 
\end{aug}

\thankstext{t1}{Supported in part by the
Natural Sciences and Engineering Research Council of Canada.}

\received{\smonth{1} \syear{2009}}
\revised{\smonth{12} \syear{2009}}

%
\begin{abstract}
In this paper we introduce a ten-parameter family of L\'evy processes
for which we obtain Wiener--Hopf factors and distribution of the
supremum process in semi-explicit form. This family allows an
arbitrary behavior of small jumps and includes processes similar to the
generalized tempered stable, KoBoL and CGMY processes. Analytically it
is characterized by the property that the characteristic exponent is a
meromorphic function, expressed in terms of beta and digamma functions.
We prove that the Wiener--Hopf factors can be expressed as infinite
products over roots of a certain transcendental equation, and the
density of the supremum process can be computed as an exponentially
converging infinite series. In several special cases when the roots can
be found analytically, we are able to identify the Wiener--Hopf factors
and distribution of the supremum in closed form. In the general case we
prove that all the roots are real and simple, and we provide
localization results and asymptotic formulas which allow an efficient
numerical evaluation. We also derive a convergence acceleration
algorithm for infinite products and a simple and efficient procedure to
compute the Wiener--Hopf factors for complex values of parameters. As a
numerical example we discuss computation of the density of the supremum
process.
\end{abstract}

%

\pdfkeywords{60G51, 60E10, Levy process,
supremum process, Wiener--Hopf factorization,
meromorphic function, infinite product}

\end{frontmatter}

\section{Introduction}\label{section_introduction}

Wiener--Hopf factorization is a powerful tool in the study of various
functionals of a L\'evy process, such as
extrema of the process, first passage time and the overshoot, the last
time the extrema was achieved, etc. These results
are very important from the theoretical point of view; for example,
they can be used to prove general theorems
about short/long time behavior (see \cite{Bertoin,Doney2007,Kyprianou} and
\cite{Sato}). However, in recent years, there has
also been a growing
interest in applications of Wiener--Hopf factorization, for example,
in Insurance Mathematics and the classical ruin problem (see
\cite{Asmussen2}) and in Mathematical Finance,
where the above-mentioned functionals are being used to describe the
payoff of a contract and
the corresponding probability distribution is used to compute its price
(see \cite{Asmussen,Boyarchenko,Mordecki3} and
\cite{Schoutens} and the references therein).

Let us summarize one of the most important results from Wiener--Hopf
factorization. Assume that
$X_t$ is a one-dimensional real-valued L\'evy process started from
$X_0=0$ and defined by a triple $(\mu,\sigma,\nu)$, where
$\mu\in{\mathbb R}$ specifies the linear component, $\sigma\ge0$ is
the volatility of the Gaussian component and
$\nu(\d x)$ is the L\'evy measure satisfying $\int_{{\mathbb R}} \min
\{
1,x^2\}
\nu(\d x) < \infty$. The characteristic exponent $\Psi(z)$ is
defined by
\[
{\mathbb E} [e^{\i zX_t} ]=e^{-t\Psi(z)},\qquad z \in {\mathbb R},
\]
and the L\'evy--Khintchine representation (see \cite{Bertoin}) tells us
that $\Psi(z)$ can be expressed in terms of the generating triple
$(\mu
,\sigma,\nu)$ as follows:
%
%
\begin{equation}\label{levy_khintchine}
\Psi(z)=\frac12 \sigma^2z^2 -\i\mu z-\int _{{\mathbb R}}
\bigl(
e^{\i zx}-1-izh(x) \bigr) \nu(\d x).
\end{equation}
Here $h(x)$ is the cut-off function, which in general can be taken to
be equal to $x\mathbf{I}_{\{|x|<1\}}$; however, in this paper we will use
$h(x)\equiv0$ (Section \ref{section_comp_poisson}) or $h(x)\equiv x$
(Sections \ref{section_results_nu3} and \ref{section_results_nu4}).

We define extrema processes
\[
S_t=\sup\{X_s \dvtx 0\le s \le t\}, \qquad    I_t=\inf\{X_s \dvtx 0\le s \le t\}
\]
introduce an exponential random variable $\tau=\tau(q)$ with parameter
$q>0$, which is independent of the process $X_t$,
and use the following notation for characteristic functions of $S_{\tau
}$ and $I_{\tau}$:
\[
\phi_q^{+}(z)={\mathbb E} [ e^{\i z S_{\tau(q)}}  ],\qquad
\phi
_q^{-}(z)={\mathbb E} [ e^{\i z I_{\tau(q)}}  ].
\]

The Wiener--Hopf factorization states that the random variables
$S_{\tau
}$ and $X_{\tau}-S_{\tau}$ are independent,
random variables $I_{\tau}$ and $X_{\tau}-S_{\tau}$ have the same
distribution; thus for $z\in{\mathbb R}$ we have
%
%
\begin{eqnarray}\label{eq_WH_factorization}
\frac{q}{q+\Psi(z)}&=&{\mathbb E} [ e^{\i zX_{\tau}}
]\nonumber\\[-8pt]\\[-8pt]
&=&{\mathbb E} [ e^{\i z S_{\tau}}  ]{\mathbb E} \bigl[
e^{\i z(X_{\tau
}-S_{\tau
})}  \bigr]=\phi_q^{+}(z)\phi_q^{-}(z).\nonumber
\end{eqnarray}
Moreover, random variable $S_{\tau}$ ($I_{\tau}$) is infinitely
divisible, positive (negative) and has no linear component
in the L\'evy--Khintchine representation (\ref{levy_khintchine}).
There also exist several integral representations for $\phi
_q^{\pm}$ in terms of ${\mathbb P}(X_t \in\d x)$ (see
\cite{Bertoin,Doney2007,Kyprianou} and \cite{Sato}) or
in terms of $\Psi(z)$ \cite{Baxter1957,Mordecki}.

The integral expressions for the Wiener--Hopf factors $\phi
_q^{\pm}$ are quite complicated; however, in the case of stable process
it is possible to obtain explicit formulas for a dense class of
parameters (see \cite{Doney1987}). It is
remarkable that in some cases
we can compute Wiener--Hopf factors explicitly with the help of
factorization identity~(\ref{eq_WH_factorization}). As an example, let us
consider the case when the L\'evy measure is of phase-type. Phase-type
distribution (see \cite{Asmussen2})
can be defined as the distribution of the first passage time of a
finite state continuous
time Markov chain. A L\'evy process $X_t$ whose jumps are phase-type
distributed enjoys the following analytical property:
its characteristic function $\Psi(z)$ is a rational function. Thus function
$q(q+\Psi(z))^{-1}$ is also a rational function, and therefore it has a
finite number of zeros/poles in the complex plane ${\mathbb C}$.
And here is the main idea: since the random variable $S_{\tau}$
($I_{\tau}$) is positive (negative) and infinitely divisible, its
characteristic function must be analytic and have no zeros in ${\mathbb C}
^{+
}$ (${\mathbb C}^{-}$), where
\[
{\mathbb C}^{+}=\{ z\in{\mathbb C}  \dvtx   \operatorname{Im}(z) > 0 \}
,\qquad
{\mathbb C}^{-}=\{ z\in{\mathbb C}  \dvtx   \operatorname{Im}(z) < 0 \}
, \qquad   \bar{\mathbb C}
^{\pm
}={\mathbb C}^{\pm}\cup{\mathbb R}.
\]
Thus we can \textit{uniquely} identify $\phi_q^{+}(z)$ [$\phi
_q^{-
}(z)$] as a rational function, which has value one at $z=0$ and whose
poles/zeros coincide with poles/zeros of $q(q+\Psi(z))^{-1}$ in
${\mathbb C}^{-}$ (${\mathbb C}^{+}$).

While L\'evy processes with phase-type jumps are very convenient
objects to work with and one can implement efficient numerical schemes,
there are some unresolved difficulties. One of them is that by
definition phase-type distribution has a smooth density on $[0,\infty)$;
in particular the density of the L\'evy measure cannot have a
singularity at zero. This means that if we want to work with a process
with infinite activity of jumps, we have to approximate its L\'evy measure
by a sequence of phase-type measures, but then the degree of rational
function $\Psi(z)$ would go to infinity and the above algorithm
for computing Wiener--Hopf factors would quickly become unfeasible.

In this paper we address this problem and discuss Wiener--Hopf
factorization for processes whose L\'evy measure can have a
singularity of arbitrary order
at zero. The main idea is quite simple: if characteristic exponent
$\Psi
(z)$ is \textit{meromorphic} in ${\mathbb C}$ and if we have sufficient
information about zeros/poles of $q+\Psi(z)$, we can still use
factorization identity (\ref{eq_WH_factorization}) essentially
in the same way as in the case of phase-type distributed jumps, except
that all the finite products will be replaced by infinite products,
and we have to be careful with the convergence issues. The main
analytical tools will be asymptotic expansion of solutions to $q+\Psi
(z)=0$ and
asymptotic results for infinite products.

The paper is organized as follows: in Section \ref
{section_comp_poisson} we introduce a simple example of a compound
Poisson process,
whose L\'evy measure has a density given by $\nu(x)=e^{\alpha x}{\mbox
{sech(x)}}$. We obtain
closed form expressions for the Wiener--Hopf factors and density of
$S_{\tau}$. Also, in this simple case we introduce many
ideas and tools which will be used in other sections. In Section \ref
{section_results_nu3} we introduce
a L\'evy process $X_t$
with jumps of infinite variation and the density of the L\'evy measure
$\nu(x)=e^{\alpha x}{\sinh(x/2)}^{-2}$. This process is a member of the
general $\beta$-family defined later in
Section \ref{section_results_nu4}; however, it is quite unique because
its characteristic
exponent $\Psi(z)$ is expressed in terms of simpler functions, and
thus all the formulas are easier and stronger results can be proved. In
this section we derive the
localization results and asymptotic expansion for the solutions of
$q+\Psi(iz)=0$, prove that
all of them are real and simple, obtain explicit formulas for sums of
inverse powers of
these solutions and finally obtain semi-explicit formulas for
Wiener--Hopf factors and distribution of supremum $S_{\tau}$.
In Section \ref{section_results_nu4} we define the ten-parameter
$\beta
$-family of L\'evy processes and derive
formulas for characteristic exponent and prove results similar to the
ones in Section \ref{section_results_nu3}.
Section \ref{section_implementation} deals with numerical issues: we
discuss acceleration of convergence of
infinite products and introduce an efficient method to compute roots of
$q+\Psi(z)$ for $q$ complex. As an example we compute
the distribution of the supremum process $S_t$.

\section{A compound Poisson process}\label{section_comp_poisson}

In this section we study a compound Poisson process $X_t$, defined by a
L\'evy measure having density
\[
\nu(x)=\frac{e^{\alpha x}}{\cosh(x)}.
\]
We take the cut-off function $h(x)$ in (\ref{levy_khintchine}) to be
equal to zero, and thus the characteristic exponent of $X_t$ is given by
%
%
\begin{equation}\label{eq_Psi_nu1}
\Psi(z)=-\int _{\mathbb R}  (e^{\i xz}-1 )\nu(x)\,\d x=
\frac{\pi}{\cos ({\pi}/2 \alpha )}-\frac{\pi
}{\cosh
({\pi}/2(z-\i\alpha) )},
\end{equation}
and the above integral can be computed with the help of formula 3.981.3
in \cite{Jeffrey2007}.
Our main result in this section is the following theorem, which
provides closed-form expressions for the Wiener--Hopf factors
and the distribution of $S_{\tau}$.
\begin{theorem}\label{thm_nu1_1} Assume that $q>0$. Define
%
%
\begin{eqnarray}\label{def_eta}
\eta&=&\frac{2}{\pi} \arccos \biggl(\frac{\pi}{q+\pi\sec
({\pi}/2 \alpha ) } \biggr), \nonumber\\[-8pt]\\[-8pt]
p_0&=&\frac{\Gamma (1/4(1-\alpha)  )\Gamma
(1/4(3-\alpha)  )}
{\Gamma (1/4(\eta-\alpha)  )\Gamma (1/4(4-\eta
-\alpha)  )}.\nonumber
\end{eqnarray}
Then for $\operatorname{Im}(z)>(\alpha-\eta)$ we have
%
%
\begin{eqnarray}\label{eq_Mtauq}
\phi_q^{+}(z)=
p_0\frac{\Gamma (1/4(\eta-\alpha-\i z)  )\Gamma
(1/4(4-\eta-\alpha-\i z)  )}
{\Gamma (1/4(1-\alpha-\i z)  )\Gamma (1/4(3-\alpha
-\i z)  )}.
\end{eqnarray}
We have ${\mathbb P}(S_{\tau}=0)=p_0$, and the density of $S_{\tau}$
is given by
%
%
\begin{eqnarray}\label{eq_density_M1}
&&\frac{\d}{\d x}{\mathbb P}(S_{\tau}\le x)\nonumber\\
&&\qquad=\frac{2p_0}{\pi}\cot
\biggl(\frac{\pi
\eta}2  \biggr) \nonumber\\
&&\qquad\quad{}\times \biggl[ \frac{\Gamma (1/4(1+\eta)  )\Gamma
(1/4(3+\eta)  )}{\Gamma (1/2 \eta )}\nonumber\\[-8pt]\\[-8pt]
&&\qquad\quad\hspace*{15.16pt}{}\times
e^{(\alpha-\eta)x}{}_2F_1 \biggl( \frac{1+\eta}4,\frac{3+\eta
}4;\frac
{\eta}2;e^{-4x}  \biggr)\nonumber\\
&&\qquad\quad\hspace*{15.16pt}{}-\frac{\Gamma (1/4(5-\eta)  )\Gamma (
1/4(7-\eta)
 )}{\Gamma (1/2(4-\eta)  )}\nonumber\\
&&\qquad\quad\hspace*{26.28pt}{}\times
e^{(\alpha-4+\eta)x}{}_2F_1 \biggl( \frac{5-\eta}4,\frac{7-\eta
}4;\frac
{4-\eta}2;e^{-4x}  \biggr)  \biggr],\nonumber
\end{eqnarray}
where ${}_2F_1(a,b;c;z)$ is the Gauss hypergeometric function.
If $q=0$ and $\alpha<0$, equation (\ref{def_eta}) implies $\eta=|
\alpha|$, and formulas (\ref{eq_Mtauq}) and
(\ref{eq_density_M1}) are still valid. In this case the random variable
$S_{\tau(0)}$ should be interpreted as $S_{\infty}=\sup\{X_s \dvtx s\ge
0 \}$.
\end{theorem}

First we will state and prove the following lemma, which will be used
repeatedly in this paper. It is a variant of the Wiener--Hopf
argument, which we have borrowed from the proof of Lemma 45.6 in
\cite{Sato}.
\begin{lemma}\label{Lemma1}
Assume we have two functions $f^{+}(z)$ and $f^{-}(z)$, such
that $f^{\pm}(0)=1$,
$f^{\pm}(z)$ are analytic in ${\mathbb C}^{\pm}$, continuous and
have no roots in $\bar{\mathbb C}^{\pm}$ and $z^{-1}\ln(f^{\pm
}(z))\to0$ as $z\to\infty$, $z\in\bar{\mathbb C}^{\pm}$.
If
%
%
\begin{equation}\label{WH_factorn_f}
\frac{q}{q+\Psi(z)}=f^{+}(z)f^{-}(z), \qquad  z\in{\mathbb R},
\end{equation}
then $f^{\pm}(z)\equiv\phi_q^{\pm}(z)$.
\end{lemma}
\begin{pf}
We define function $F(z)$ as
\[
F(z)=\cases{
\dfrac{\phi_q^{-}(z)}{f^{-}(z)}, &\quad if
$z \in\bar{\mathbb C}^{-}$, \cr
\dfrac{f^{+}(z)}{\phi_q^{+}(z)}, &\quad if
$z \in\bar{\mathbb C}^{+}$.}
\]
Function $F(z)$ is well defined for $z$ real due to (\ref
{WH_factorn_f}) and (\ref{eq_WH_factorization}).
Using properties of $\phi_q^{\pm}$ and $f^{\pm}$ we
conclude that $F(z)$ is analytic in ${\mathbb C}^{+}$ and
${\mathbb
C}^{-}$ and continuous in ${\mathbb C}$, and therefore by analytic
continuation (see Theorem 16.8 on page 323 in \cite{Rudin1986}) it must
be analytic in the entire complex plane.
Moreover, by construction function $F(z)$ has no zeros in ${\mathbb
C}$, thus
its logarithm is also an entire function. All that is left to do is to
prove that function $\ln(F(z))$ is constant.

Using integration by parts and formula (\ref{levy_khintchine}) one
could prove the following result:
if $\xi$ is an infinitely divisible positive random variable with no
drift and $\Psi_{\xi}(z)$ is its characteristic exponent,
then $z^{-1} \Psi_{\xi}(z)\to0$ as $z\to\infty$, $z\in\bar
{\mathbb C}
^{+
}$ (this statement is similar to Proposition 2 in \cite{Bertoin}). Thus
\[
z^{-1} \ln(\phi_q^{\pm}(z)) \to0,  \qquad   z \to\infty,\qquad
z\in\bar{\mathbb C}^{\pm}.
\]
Since functions $f^{\pm}$ also satisfy the above conditions, we
find that $z^{-1} \ln(F(z))\to0$ as $| z |\to\infty$ in the
entire complex plane. Thus we have an analytic function $\ln(F(z))$
which grows slower than $|z|$ as $z\to\infty$, and therefore we can
conclude that this function
must be constant (a rigorous way to prove this is to apply Cauchy's
estimates, see Proposition 2.14 on page 73 in \cite{Conway1978}). The
value of this constant is easily seen to be zero, since $f^{\pm
}(0)=\phi_q^{\pm}(0)=1$.
\end{pf}
\begin{pf*}{Proof of Theorem \ref{thm_nu1_1}}
Using expression (\ref{eq_Psi_nu1}) for $\Psi(z)$ we find that function
$q(q+\Psi(z))^{-1}$ has simple zeros at $\{\i(1+\alpha+4n),\i
(3+\alpha
+4n)\}$ and simple poles at $\{\i(\alpha+\eta+4n),\i(\alpha-\eta
+4n)\}
$, where $n\in{\mathbb Z}$ and $\eta$ is defined by (\ref{def_eta}).
Next we check that $|\alpha|< \eta<1$ and
define function $f^{+}(z)$ as product over all zeros/poles lying in
${\mathbb C}^{-}$
%
%
\begin{equation}\label{def_f_plus}
f^{+}(z)=\prod _{n\ge0} \frac
{ (1-{\i z}/({4n+1-\alpha})) (1-{\i
z}/({4n+3-\alpha
}) )}
{ (1-{\i z}/({4n+\eta-\alpha} )) (1-{\i
z}/({4n+4-\eta
-\alpha}) )}
\end{equation}
and similarly $f^{-}(z)$ as product over zeros/poles in ${\mathbb
C}^{+}$. It is easy to see that the product
converges uniformly on compact subsets of ${\mathbb C}\setminus\i
{\mathbb R}$ since
each term is $1+O(n^{-2})$
(see Corollary 5.6 on page 166 in \cite{Conway1978} for sufficient
conditions for the absolute convergence of infinite products).
The fact that $f^{+}(z)$ is equal to the right-hand side of formula
(\ref{eq_Mtauq}) can be seen by applying the following result
from \cite{Erdelyi1955V3}:
%
%
\begin{equation}
\prod _{n \ge0} \frac
{1+{x}/({n+a})}{1+{x}/({n+b})}
=\frac{\Gamma(a)\Gamma(b+x)}{\Gamma(b)\Gamma(a+x)}.
\end{equation}
The formula for $f^{-}(z)$ is identical to (\ref{eq_Mtauq}) with
$(z, \alpha)$
replaced by $(-z, -\alpha)$.

Now we will prove that $f^{\pm}(z)\equiv\phi_q^{\pm
}(z)$. First, using the reflection formula for the gamma function
(formula 8.334.3 in \cite{Jeffrey2007}), one can check that for $z \in
{\mathbb R}$ functions $f^{\pm}(z)$ satisfy factorization identity
(\ref
{WH_factorn_f}).
Next, using the following asymptotic expression (formula 6.1.47 in
\cite{AbramowitzStegun}):
%
%
\begin{equation}\label{eq_gamma_asympt}
\frac{\Gamma(a+x)}{\Gamma(b+x)}=x^{a-b}+O(x^{a-b-1}),
\end{equation}
we conclude that $z^{-1} \ln(f^{\pm}(z)) \to0$ as $z \to
\infty
$, $z\in\bar{\mathbb C}^{\pm}$, and thus all the conditions of
Lemma \ref{Lemma1} are satisfied, and we conclude that $f^{\pm
}(z)\equiv\phi_q^{\pm}(z)$.

To derive formula (\ref{eq_density_M1}) for the density of $S_{\tau}$
we use equations (\ref{eq_Mtauq}) and (\ref{eq_gamma_asympt})
to find that ${\mathbb E} [e^{-\zeta S_{\tau}}  ]=\phi_q^{+
}(\i
\zeta)\to p_0$ as $\zeta\to\infty$, where $p_0$ is given by (\ref
{def_eta}). This implies
that distribution of $S_{\tau}$ has an atom at $x=0$ (which should not
be surprising since $X_t$ is a compound Poisson process), and
${\mathbb P}(S_{\tau}=0)=p_0$. The density of $S_{\tau}$ can be
computed by the
inverse Fourier transform
\begin{eqnarray*}
&&\frac{\d}{\d x}{\mathbb P}(S_{\tau}\le x)\\
&&\qquad=\frac{1}{2\pi} \int
_{\mathbb
R}  [\phi_q^{+}(z)-p_0 ] e^{-\i xz}\,\d z\\
&&\qquad=
\frac{p_0}{2\pi} \int _{\mathbb R}
\biggl[\frac{\Gamma (1/4(\eta-\alpha-\i z)  )\Gamma
(1/4(4-\eta-\alpha-\i z)  )}
{\Gamma (1/4(1-\alpha-\i z)  )\Gamma (
1/4(3-\alpha
-\i z)  )}-1 \biggr]e^{-\i xz} \,\d z.
\end{eqnarray*}
Formula (\ref{eq_density_M1}) is obtained from the above expression by
replacing the contour of integration by $\i c+{\mathbb R}$, letting
$c\to-\infty$ and evaluating the residues at
$z \in\{-\i(4n+\eta-\alpha), -\i(4n+4-\eta-\alpha)\}$ for $n\ge0$.
Evaluating the residues can be made easier by using the reflection
formula for the gamma function.
\end{pf*}
\begin{remark}
There are other examples of L\'evy measures $\nu(x)\,dx$, which have
finite total mass (and thus can define a process with
a finite intensity of jumps), and for which the characteristic exponent
is a
simple meromorphic function. These are two examples based on theta
functions (see Section 8.18 in \cite{Jeffrey2007} for definition and
properties of theta functions):
\begin{eqnarray*}
\nu_1(x)&=&e^{-\alpha x} \theta_2 (0 ,e^{-x}  )=e^{-\alpha
x} \biggl[2\sum _{n\ge0} e^{-(n+1/2)^2x} \biggr], \\
\nu_2(x)&=& e^{-\alpha x} \theta_3 (0 ,e^{-x}
)=e^{-\alpha
x} \biggl[1+2\sum _{n\ge0} e^{-n^2x} \biggr].
\end{eqnarray*}
These two jump densities are defined on $x>0$, they decay exponentially
as $x\to+\infty$ and behave as $x^{-1/2}$ as $x\to0^+$;
thus the total mass is finite. The Fourier transform of these functions
can be computed using formulas 6.162 in \cite{Jeffrey2007}
\begin{eqnarray*}
\int _{0}^{\infty} e^{\i xz} \nu_1(x)\,\d x&=&\frac{\pi
}{\sqrt
{\alpha-\i z}} \tanh \bigl(\pi\sqrt{\alpha-\i z} \bigr), \\
\int _{0}^{\infty} e^{\i xz} \nu_2(x)\,\d x&=&\frac{\pi
}{\sqrt
{\alpha-\i z}} \coth \bigl(\pi\sqrt{\alpha-\i z} \bigr).
\end{eqnarray*}
Unfortunately equation $q+\Psi(z)=0$ cannot be solved explicitly which
implies that we cannot obtain closed form results
as in Theorem \ref{thm_nu1_1}; however, these processes could be
treated using methods presented in the next sections.
\end{remark}

\section{A process with jumps of infinite variation}\label{section_results_nu3}

In this section we study a L\'evy process $X_t$, defined by a triple
$(\mu,\sigma,\nu)$, where the density of the L\'evy measure is given by
\[
\nu(x)=\frac{e^{\alpha x}}{ [\sinh({x}/2) ]^2}
\]
with $|\alpha|<1$ (it is a L\'evy measure of a Lamperti-stable
process with characteristics $(1,1+\alpha,1-\alpha)$, see \cite{KyPaRi}).
The jump part of $X_t$ is similar to the normal inverse Gaussian
process (see \cite{Barndorff1997,Cont}), as it is also a process of infinite
variation,
the jump measure decays exponentially as $|x| \to\infty$ and has a
$O(x^{-2})$ singularity at $x=0$. Note that since
the L\'evy measure has exponential tails we can take the cut-off
function $h(x)\equiv x$ in (\ref{levy_khintchine}).

By definition process $X_t$ has three parameters. However, if we want to
achieve greater generality for modeling purposes, we could introduce
two additional scaling parameters $a$ and $b>0$ and define
a process $Y_t=aX_{bt}$, thus obtaining a five parameter family of L\'
evy processes.
\begin{proposition}
The characteristic exponent of $X_t$ is given by
%
%
\begin{equation}\label{eq_psi_nu3}
\Psi(z)=\tfrac12 \sigma^2z^2+i\rho z
+4\pi(z-i\alpha)\coth \bigl(\pi(z-i\alpha) \bigr)-4\gamma,
\end{equation}
where
\[
\gamma=\pi\alpha\cot (\pi\alpha ),  \qquad   \rho=4\pi
^2 \alpha
+\frac{4\gamma(\gamma-1)}{\alpha}-\mu.
\]
\end{proposition}
\begin{pf}
We start with the series representation valid for $x>0$,
%
%
\begin{equation}\label{sinh_series}
\biggl[\sinh \biggl(\frac{x}2 \biggr) \biggr]^{-2}=4\frac
{e^{-x}}{
(1-e^{-x} )^2}=4\sum _{n\ge1} ne^{-nx},
\end{equation}
which can be easily obtained
using binomial series or by taking derivative of a geometric series.
The infinite series in (\ref{sinh_series}) converges uniformly on
$(\varepsilon,\infty)$
for every $\varepsilon>0$, thus
\begin{eqnarray*}
&&
\int _{0}^{\infty}  ( e^{\i zx}-1-\i zx ) \frac
{e^{\alpha x}}{\sinh({x}/2)^2}\,\d x\\
&&\qquad=
4\sum _{n\ge1}  \biggl[ \frac{n}{n-\alpha-\i z}-\frac
{n}{n-\alpha
}- \frac{\i nz}{(n-\alpha)^2} \biggr]\\
&&\qquad=
4\sum _{n\ge1}  \biggl[\frac{\alpha+\i z}{n-\alpha-\i
z}-\frac
{\alpha+\i z}{n-\alpha}- \frac{\i\alpha z}{(n-\alpha)^2}  \biggr].
\end{eqnarray*}
The integral in the L\'evy--Khintchine representation (\ref
{levy_khintchine}) for $\Psi(z)$ can now be computed as
\begin{eqnarray*}
&&
\int _{0}^{\infty}  ( e^{\i zx}-1-\i zx ) \frac
{e^{\alpha x}}{\sinh(x)^2}\,\d x+
\int _{0}^{\infty}  ( e^{-\i zx}-1+\i zx ) \frac
{e^{-\alpha x}}{\sinh(x)^2}\,\d x\\
&&\qquad=
8(\alpha+\i z)^2\sum _{n \ge1}\frac{1}{n^2-(\alpha+\i z)^2}-
8(\alpha+\i z)\alpha\sum _{n \ge1}\frac{1}{n^2-\alpha^2}\\
&&\qquad\quad{}  -
4\i\alpha z \sum _{n \in{\mathbb Z}} \frac{1}{(n-\alpha
)^2}+\frac
{4\i z}{\alpha}.
\end{eqnarray*}
To complete the proof we need to use the following well-known series
expansions (see formulas 1.421.4 and 1.422.4 in \cite{Jeffrey2007})
\begin{eqnarray*}
\coth(\pi x)&=&\frac{1}{\pi x}+\frac{2x}{\pi}\sum _{n \ge
1}\frac
{1}{x^2+n^2}, \\
1+\cot(\pi x)^2&=&\operatorname{cosec}(\pi x)^2=\frac{1}{\pi^2} \sum
_{n \in{\mathbb Z}} \frac{1}{(n-x)^2}.
\end{eqnarray*}
\upqed\end{pf}

Note that it is impossible to find solutions to $q+\Psi(z)=0$
explicitly in the general case, even though the characteristic exponent
$\Psi(z)$ is quite simple. It is remarkable that in some special cases,
when $\sigma=0$ and parameters $\mu$, $\alpha$ and $q$ satisfy certain
conditions, we can still obtain closed-form results.
Below we present just one example of this type.
\begin{proposition}\label{prop_explicit_nu3} Assume that $\sigma
=\alpha
=0$. Define
%
%
\begin{equation}
\eta=\frac{1}{\pi}\operatorname{arccot} \biggl( \frac{\mu}{4\pi
} \biggr).
\end{equation}
Then Wiener--Hopf factor $\phi_q^{+}(z)$ can be computed in closed
form when $q=4$,
%
%
\begin{equation}
\phi_4^{+}(z)=\frac{\Gamma(\eta-\i z)}{\Gamma(\eta)\Gamma
(1-\i z)}.
\end{equation}
The density of $S_{\tau(4)}$ is given by
\[
\frac{\d}{\d x}   {\mathbb P}\bigl(S_{\tau(4)}\le x\bigr) =\frac{\sin(\pi
\eta
)}{\pi}
 ( e^x-1  )^{-\eta}.
\]
\end{proposition}

The proof of Proposition \ref{prop_explicit_nu3} is identical to the
proof of Theorem \ref{thm_nu1_1}.
\begin{remark}
This result is very similar to Proposition 1 in \cite{CaCha}, where the
authors are able to compute the law
of $I_{\tau(q)}$ in closed form only for a single value of~$q$, and
this law is essentially identical to the distribution of $S_{\tau(4)}$.
This coincidence seems to be rather surprising, since these
propositions study different processes: our Proposition \ref{prop_explicit_nu3}
is concerned with a Lamperti-stable process having characteristics
$(1,1,1)$ (see \cite{KyPaRi}) and completely arbitrary drift, while
Proposition~1 in \cite{CaCha} studies a Lamperti-stable process with
characteristics $(\alpha,1,\alpha)$ but with no freedom
in specifying the drift, which must be uniquely expressed in terms of
parameters of the L\'evy measure.
\end{remark}

The following theorem is one of the main results in this section. It
describes various properties of solutions to equation $q+\Psi(z)=0$,
which will be used later to compute Wiener--Hopf factors and the
distribution of the supremum process.
\begin{theorem}\label{thm_nu3_1} Assume that $q>0$ and that $\Psi(z)$
is given by (\ref{eq_psi_nu3}).

\begin{longlist}
\item Equation $q+\Psi(i\zeta)=0$ has infinitely many solutions,
all of which are real and simple.
They are located as follows:
%
%
\begin{eqnarray}\label{eq_localization}
\zeta_0^{-} &\in& (\alpha-1,0),\nonumber\\
\zeta_0^{+} &\in& (0,\alpha+1), \nonumber\\[-8pt]\\[-8pt]
\zeta_n &\in& (n+\alpha, n+\alpha+1),  \qquad   n\ge1 \nonumber\\
\zeta_n &\in& (n+\alpha-1, n+\alpha),  \qquad   n\le-1.\nonumber
\end{eqnarray}
\item
If $\sigma\ne0$ we have as $n\to\pm\infty$
%
%
\begin{eqnarray}\label{eq_asympt_roots1_1}
\zeta_n&=&(n+\alpha)+\frac{8}{\sigma^2}(n+\alpha)^{-1}\nonumber\\[-8pt]\\[-8pt]
&&{}-\frac{8}{\sigma^2} \biggl(\frac{2\rho}{\sigma^2}+\alpha
\biggr)(n+\alpha)^{-2}+O(n^{-3}).\nonumber
\end{eqnarray}
\item
If $\sigma=0$ we have as $n\to\pm\infty$
%
%
\begin{eqnarray}\label{eq_asympt_roots1_2}
\zeta_{n+\delta}&=&(n+\alpha+\omega_0)+c_0(n+\alpha+\omega
_0)^{-1}\nonumber\\[-8pt]\\[-8pt]
&&{}-\frac{c_0}{\rho} (4\gamma-q-4\pi^2 c_0 )
(n+\alpha+\omega_0)^{-2}+O(n^{-3}),\nonumber
\end{eqnarray}
where
\[
c_0=-\frac{4(4\gamma-q+\alpha\rho) }{16\pi^2+\rho^2},\qquad
\omega
_0=\frac{1}{\pi}\operatorname{arccot} \biggl(\frac{\rho}{4\pi} \biggr)
\]
and $\delta\in\{-1,0,1\}$ depending on the signs of $n$ and $\rho$.
\item Function $q(q+\Psi(z))^{-1}$ can be factorized as follows:
%
%
\begin{equation}\label{eq_big_factorization}
\frac{q}{q+\Psi(z)}=\frac{1}{ (1+{\i z}/{\zeta_0^{+
}}
) (1+{\i z}/{\zeta_0^{-}}  )}
\prod _{| n |\ge1} \frac{1+{\i z}/({n+\alpha
})}{1+
{\i z}/{\zeta_n}},
\end{equation}
where the infinite product converges uniformly on the compact subsets
of the complex plane excluding zeros/poles of $q+\Psi(z)$.
\end{longlist}
\end{theorem}

First we need to prove the following technical result.
\begin{lemma}\label{lemmma_asymp_for_product}
Assume that $\alpha$ and $\beta$ are not equal to a negative integer,
and $b_n=O(n^{-\varepsilon_1})$ for some $\varepsilon_1>0$ as $n\to\infty
$. Then
%
%
\begin{equation}\label{eq_inf_prod}
\prod _{n \ge0} \frac{1+{z}/({n+\alpha})}{1+
{z}/({n+\beta
+b_n})} \approx C z^{\beta-\alpha}
\end{equation}
as $z\to\infty$, $|{\operatorname{arg}}(z) |< \pi-\varepsilon_2<\pi
$, where
$C=\frac{\Gamma(\alpha)}{\Gamma(\beta)}\prod_{n \ge0}
(1+\frac
{b_n}{n+\beta} )$.
\end{lemma}
\begin{pf}
First we have to justify absolute convergence of infinite products.
A~product in the left-hand side of (\ref{eq_inf_prod}) converges since
each term is
$1+O(n^{-2})$ and the infinite product in the definition of constant
$C$ converges since each term is $1+O(n^{-1-\varepsilon_1})$
(see Corollary 5.6 on page 166 in \cite{Conway1978} for sufficient
conditions for the absolute convergence of infinite products). Thus we
can rewrite the left-hand side of (\ref{eq_inf_prod}) as
%
%
\begin{eqnarray}\label{asympt_proof1}\qquad
\prod _{n \ge0} \frac{1+{z}/({n+\alpha})}{1+
{z}/({n+\beta+b_n})}&=&
\prod _{n \ge0} \frac{1+{z}/({n+\alpha})}{1+
{z}/({n+\beta})}
\prod _{n \ge0} \frac{1+{z}/({n+\beta})}{1+
{z}/({n+\beta
+b_n})}\nonumber\\[-8pt]\\[-8pt]
&=&
C\frac{\Gamma(\beta+z)}{\Gamma(\alpha+z)} \prod _{n\ge0}
\frac
{z+n+\beta}
{z+n+\beta+b_n}.\nonumber
\end{eqnarray}
The ratio of gamma functions gives us the leading asymptotic term
$z^{\beta-\alpha}$ due to~(\ref{eq_gamma_asympt}). Now we need to prove that the last infinite
product in (\ref{asympt_proof1}) converges to one
as $z\to\infty$, $|{\operatorname{arg}}(z) |< \pi-\varepsilon_2<\pi$.
We take the logarithm of this product and estimate it as
\begin{eqnarray*}
\biggl| \sum _{n\ge1} \ln \biggl(\frac{z+n+\beta}
{z+n+\beta+b_n}  \biggr)  \biggr| &=&
\biggl| \sum _{n\ge1} \ln \biggl(1+\frac{b_n}
{z+n+\beta}  \biggr)  \biggr| \nonumber\\
&\le&
\sum _{n\ge1} \ln \biggl(1+\frac{| b_n |}
{| z+n+\beta|}  \biggr) \nonumber\\
&\le& \sum _{n\ge1} \frac{| b_n |}
{| z+n+\beta|}\\
&\le& A\sum _{n\ge1} \frac{1}
{n^{\varepsilon_1}| z+n+\beta|},
\end{eqnarray*}
where we have used the fact that $\ln(1+x)<x$ for $x>0$ and $| b_n
|< An^{-\varepsilon_1}$ for some $A>0$.
Since $|{\operatorname{arg}}(z) |< \pi-\varepsilon_2<\pi$ we have for
$z$ sufficiently large
$|z + n + \beta| > \max\{1,|n- | z+\beta| |\}$. Let $m=[|z+\beta|]$,
where $[x]$ denotes the integer part of $x$.
Then
%
%
\begin{equation}\label{two_series}
\sum _{n\ge1} \frac{1}
{n^{\varepsilon_1}|z+n+\beta|}<\sum _{n=1}^{m} \frac{1}
{n^{\varepsilon_1}(m+1-n)}+\sum _{n=m+1}^{\infty} \frac{1}
{n^{\varepsilon_1}(n-m)}.
\end{equation}

The first series in the right-hand side of (\ref{two_series}) converges
to zero as $m\to\infty$, since
\begin{eqnarray*}
\sum _{n=1}^{m} \frac{1}
{n^{\varepsilon_1}(m+1-n)}&=&\sum _{n=1}^{[\sqrt{m}]} \frac{1}
{n^{\varepsilon_1}(m+1-n)}+\sum _{n=[\sqrt{m}]+1}^m \frac{1}
{n^{\varepsilon_1}(m+1-n)} \\
&<&\frac{[\sqrt{m}]}{m+1-[\sqrt{m}]}+m^{-\varepsilon_1/2} \sum
_{n=[\sqrt{m}]+1}^m \frac{1}
{(m+1-n)} \\ &<& \frac{[\sqrt{m}]}{m+1-[\sqrt{m}]}+m^{-\varepsilon
_1/2} \ln(m).
\end{eqnarray*}

The second series in the right-hand side of (\ref{two_series}) can be
rewritten as $\sum_{n=1}^{\infty}
(n+m)^{-\varepsilon_1}n^{-1}$, and we see that it is a convergent series
of positive terms, where each term converges to zero as $m \to\infty$.
By considering its partial sums it is easy to prove that the series
itself must converge to zero as $m \to\infty$.
\end{pf}
\begin{pf*}{Proof of Theorem \ref{thm_nu3_1}}
The proof consists of three steps. The first step is to study solutions
to equation $q+\Psi(\i\zeta)=0$.
We will produce a sequence of ``obvious'' solutions $\zeta_n$ and study
their asymptotics as $n\to\pm\infty$.
Note that this first step requires quite demanding computations, which
can be made much more enjoyable if one uses a symbolic computation package.
The second step is to represent the function $q(q+\Psi(z))^{-1}$ as a
general infinite product, which includes poles of $\Psi(z)$, zeros of
$q+\Psi(z)$ (given by ``obvious'' ones
$\{\i\zeta_0^{\pm},\i\zeta_n\}$
and possibly some ``unaccounted'' zeros) and an exponential factor. Our
main tool will be Hadamard theorem (Theorem 1, page 26, in
\cite{Levin1996} or Theorem 3.4, page 289, in \cite{Conway1978}). We produce
entire functions $P(z)$ and $Q(z)$, such that
$P(z)$ has zeros at poles of $\Psi(z)$ and $q(q+\Psi
(z))^{-1}=P(z)/Q(z)$. After studying the growth rate of $P(z)$ and $Q(z)$
we apply Hadamard theorem and obtain an infinite product for function
$Q(z)$ [function $P(z)$ will have an explicit infinite product]. These
results give us an infinite product for $q(q+\Psi(z))^{-1}$. The last
step is to prove the absence of exponential factor and ``unaccounted''
zeros in this infinite product, and here the main tool will be
asymptotic relation (\ref{eq_inf_prod}) for infinite products provided
by Lemma \ref{lemmma_asymp_for_product}.

First we will prove localization result (\ref{eq_localization}). We use
(\ref{eq_psi_nu3}) to rewrite equation $q+\Psi(\i\zeta)=0$ as
%
%
\begin{equation}\label{eq_Psiq0}
4\pi(\zeta-\alpha) \cot\bigl(\pi(\zeta-\alpha)\bigr)-(\rho+\mu)\zeta
-4\gamma
=\tfrac12 \sigma^2\zeta^2-\mu\zeta-q.
\end{equation}
Note that we have separated the jump part of $\Psi(z)$ on the left-hand
side and the diffusion part on the right-hand side of (\ref
{eq_Psiq0}). See Figure \ref{fig_plot}, where the jump part is
represented by black line and diffusion part by grey dotted line.

%
%
\begin{figure}[b]

\includegraphics{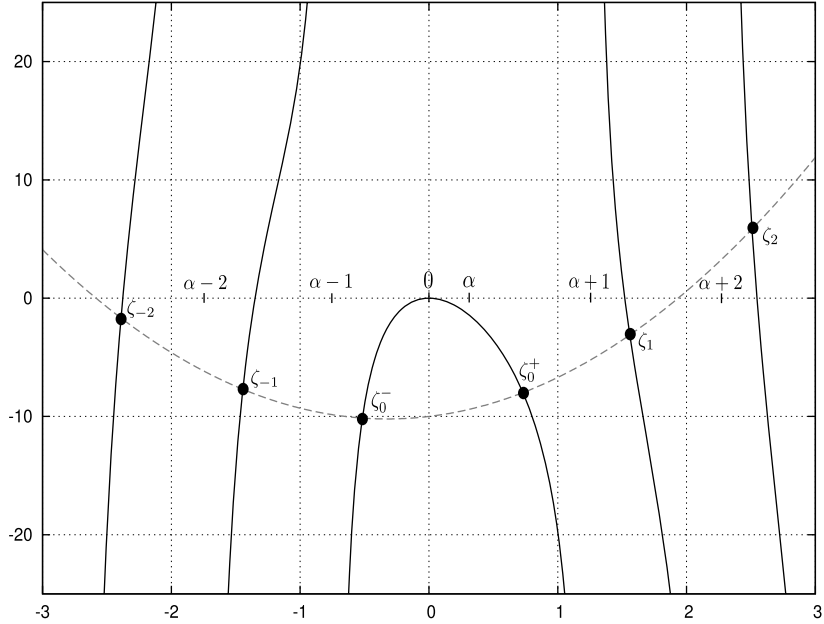}

\caption{Illustration of the proof of Theorem \protect\ref{thm_nu3_1}.}
\label{fig_plot}
\end{figure}

The left-hand side of (\ref{eq_Psiq0}) is zero at $\zeta=0$
and goes to $-\infty$ as $\zeta\nearrow\alpha+1$ or $\zeta\searrow
\alpha-1$ (see Figure \ref{fig_plot}). The right-hand side is negative
at $\zeta=0$ and continuous everywhere; thus we have at least one solution
$\zeta_0^{+} \in(0,\alpha+1)$ and at least one solution $\zeta
_0^{-} \in(\alpha-1,0)$. In fact
it is easy to prove that we have
\textit{exactly} one solution on each of these intervals, since $4\pi
(\zeta-\alpha) \cot(\pi(\zeta-\alpha))$ is a concave function on
$(\alpha-1,\alpha+1)$, while
$\frac12 \sigma^2 \zeta^2-\mu\zeta-q$ is convex.

Next, for $n\ne0$ we have
\begin{eqnarray*}
&& 4\pi(\zeta-\alpha) \cot\bigl(\pi(\zeta-\alpha)\bigr) \nearrow+\infty
\qquad\mbox{as }    \zeta\nearrow\alpha-n,
\zeta\searrow\alpha+n, \\
&& 4\pi(\zeta-\alpha) \cot\bigl(\pi(\zeta-\alpha)\bigr) \searrow-\infty
\qquad\mbox{as }    \zeta\searrow\alpha-n,
\zeta\nearrow\alpha+n,
\end{eqnarray*}
thus there must exist at least one zero $\zeta_n$ on each interval
$(n+\alpha, n+\alpha+1)$,
$(n+\alpha-1, n+\alpha)$.

Next we will prove the asymptotic expansion (\ref{eq_asympt_roots1_1}).
Since we have assumed that $\sigma\ne0$ we can rearrange the terms
in (\ref{eq_Psiq0}) to obtain
%
%
\begin{eqnarray}\label{eq_asymp1_proof1}
\frac{1}{\pi} \tan\bigl(\pi(\zeta-\alpha)\bigr)&=&\frac{4(\zeta-\alpha
)}{1/2
\sigma^2 \zeta^2+\rho\zeta+4\gamma-q}
\nonumber\\
&=&\frac{8}{\sigma^2}\zeta^{-1}
\biggl[\frac{1-\alpha\zeta^{-1}}{1+2\rho\sigma^{-2} \zeta
^{-1}+O(\zeta
^{-2})}  \biggr]\\
&=&
\frac{8}{\sigma^2}\zeta^{-1}-
\frac{8}{\sigma^2} \biggl(\frac{2}{\rho}+\alpha \biggr)\zeta
^{-2}+O(\zeta^{-3}).\nonumber
\end{eqnarray}
The main idea in the above calculation is to expand the rational
function in the Taylor series centered at $\zeta=\infty$.
Now, the right-hand side of (\ref{eq_asymp1_proof1}) is small when
$\zeta$ is large, and thus the solution to (\ref{eq_asymp1_proof1})
should be close to the solution of $\tan(\pi(\zeta-\alpha))=0$,
which implies
%
%
\begin{equation}\label{eq_zeta_omega}
\zeta= n+\alpha+ \omega
\end{equation}
and $\omega=o(1)$ as $n\to\infty$. Next we expand the right-hand side
of (\ref{eq_asymp1_proof1}) in powers of $w$ as
\[
\frac{1}{\pi} \tan\bigl(\pi(\zeta-\alpha)\bigr)=\frac{1}{\pi} \tan(\pi
\omega
)=\omega+O(\omega^3)
\]
and, using the first two terms of the Maclaurin series for $\zeta^{-1}$
in powers of $\omega$
\[
\zeta^{-1}=(n+\alpha+\omega)^{-1}=(n+\alpha)^{-1}-\omega(n+\alpha
)^{-2}+O(\omega^2n^{-3}),
\]
we are able to
rewrite (\ref{eq_asymp1_proof1}) as
\[
\omega+O(\omega^3)=\frac{8}{\sigma^2} \bigl((n+\alpha)^{-1}+\omega
(n+\alpha)^{-2} \bigr)-
\frac{8}{\sigma^2} \biggl(\frac{2}{\rho}+\alpha \biggr)(n+\alpha
)^{-2}+O(n^{-3}).
\]
Asymptotic expansion (\ref{eq_asympt_roots1_1}) follows easily from the
above formula and (\ref{eq_zeta_omega}).

If $\sigma=0$, equation (\ref{eq_asymp1_proof1}) has to be modified
as follows:
%
%
\begin{eqnarray}\label{eq_asymp1_proof2}
\frac{1}{\pi} \tan\bigl(\pi(\zeta-\alpha)\bigr)&=&\frac{4(\zeta-\alpha
)}{\rho
\zeta+4\gamma-q}=
\frac{4}{\rho}  \biggl[ \frac{1-\alpha\zeta^{-1}}{1+(4\gamma
-q)\rho
^{-1}\zeta^{-1}}  \biggr]\nonumber\\
&=&
\frac{4}{\rho}-\frac{4}{\rho^2}(4\gamma-q+\alpha\rho) \zeta^{-1}\\
&&{}+
\frac{4(4\gamma-q)}{\rho^3} (4\gamma-q+\alpha\rho) \zeta
^{-2}+O(\zeta^{-3}),\nonumber
\end{eqnarray}
where again we have expanded the rational function in the Taylor series
centered at $\zeta=\infty$.
As before, when $\zeta$ is large the solution of (\ref{eq_Psiq0})
should be close to the solution of
\[
\frac{1}{\pi} \tan\bigl(\pi(\zeta-\alpha)\bigr)=\frac{4}{\rho},
\]
and thus we should expand both sides of (\ref{eq_asymp1_proof2}) in the
Taylor series centered at the solution to the above equation. We define
$\omega$ as
%
%
\begin{equation}\label{eq_zeta_omega2}
\zeta=n+\alpha+\frac{1}{\pi}\arctan \biggl(\frac{4\pi}{\rho
} \biggr)+\omega
=n+\alpha+\omega_0+\omega,
\end{equation}
and again $\omega=o(1)$ as $n \to\infty$. To expand the left-hand side
of (\ref{eq_asymp1_proof2}) in power series in $\omega$ we use an
addition formula for $\tan(\cdot)$ and find that
%
%
\begin{eqnarray}\label{eq_tan2}\quad
\frac{1}{\pi} \tan\bigl(\pi(\zeta-\alpha)\bigr)&=&\frac{1}{\pi}
\tan
 \biggl(\arctan \biggl(\frac{4\pi}{\rho} \biggr)+\pi\omega \biggr)
\nonumber\\
&=&
\frac{1}{\pi}\frac{{4\pi}/{\rho}+\tan(\pi\omega)}{1-
{4\pi}/{\rho}\tan(\pi\omega)}
=\frac{1}{\pi}\frac{{4\pi}/{\rho}+\pi\omega+O(\omega
^3)}{1-
{4\pi}/{\rho}\pi\omega+O(\omega^3)} \nonumber\\[-8pt]\\[-8pt]
&=&\frac{4}{\rho}+
\frac{1}{\rho^2} (16 \pi^2+\rho^2 ) \omega\nonumber\\
&&{} +\frac{4 \pi
^2}{\rho
^3} (16 \pi^2+\rho^2 ) \omega^2+O(\omega^3).\nonumber
\end{eqnarray}
Again, we use (\ref{eq_zeta_omega2}) to obtain the Maclaurin series of
$\zeta^{-1}$ in powers of $\omega$
\begin{eqnarray*}
\zeta^{-1}&=&(n+\alpha+\omega_0+\omega)^{-1}\nonumber\\
&=&(n+\alpha+\omega_0)^{-1}-\omega(n+\alpha+\omega
_0)^{-2}+O(\omega^2n^{-3}).
\end{eqnarray*}
Using (\ref{eq_tan2}) and the above expansion we can rewrite
(\ref{eq_asymp1_proof2}) as
\begin{eqnarray*}
&& (16 \pi^2+\rho^2 )  \biggl[ \omega+\frac{4 \pi
^2}{\rho}\omega
^2 \biggr]\\
&&\qquad=
-4(4\gamma-q+\alpha\rho) \bigl((n+\alpha
+\omega
_0)^{-1}-\omega(n+\alpha+\omega_0)^{-2} \bigr)
\\
&&\qquad\quad{}  +\frac{4(4\gamma-q)}{\rho} (4\gamma-q+\alpha\rho)
(n+\alpha+\omega_0)^{-2}+O(n^{-3})+O(\omega^3),
\end{eqnarray*}
and from this equation we obtain the second asymptotic expansion (\ref
{eq_asympt_roots1_2}).

Now we are ready to prove the factorization identity (\ref
{eq_big_factorization}) and the fact that all the
zeros of $q+\Psi(\i\zeta)$ are real and simple and that there are no
other zeros except for the ones described in (\ref{eq_localization}).
First we need to find
an analytic function $P(z)$ such that $P(0)=1$ and which has zeros at
all poles of $\Psi(z)$ (with the same multiplicity).
The choice is rather obvious due to (\ref{eq_psi_nu3}):
%
%
\begin{equation}\label{def_P}
P(z)=\frac{\alpha}{\sin(\pi\alpha)}\times\frac{\sinh(\pi(z-\i
\alpha
))}{ z-\i\alpha}.
\end{equation}
By definition, the function
%
%
\begin{equation}\label{def_Q}
Q(z)=q^{-1}\bigl(q+\Psi(z)\bigr)P(z)
\end{equation}
is also analytic in the entire complex plane.

Next, using the definition of $P(z)$ (\ref{def_P}) and $Q(z)$ (\ref
{def_Q}) we check that $Q(z)=0$ if and only if $q+\Psi(z)=0$.
We have proved already that the zeros of $q+\Psi(\i\zeta)$ include
$\zeta_n, \zeta_0^{\pm}$; however, some of them might
have multiplicity greater than one, and there also might exist other
roots (real and/or complex). Let us denote the set of
these unaccounted roots (counting with multiplicity) as
${\mathfrak Z}$. Using asymptotic expansions given by equations
(\ref{eq_asymp1_proof1}) and (\ref{eq_asymp1_proof2}) one can easily
prove that ${\mathfrak Z}$ is a finite set (possibly empty).

Using equations (\ref{def_P}), (\ref{def_Q}) and (\ref{eq_psi_nu3}) we
obtain an explicit formula for $Q(z)$ from which it easily follows
that $Q(z)$ has order equal to one, which means that one is the least
lower bound of all $\gamma>0$ such that
$Q(z)=O(\exp(|z|^{\gamma}))$ as $z\to\infty$; the rigorous definition
can be found in \cite{Levin1996}, page 4 or Chapter 11 in \cite
{Conway1978}. Since $Q(z)$ has order equal to one, we can use the
Hadamard theorem (see Theorem 1, page~26, in \cite{Levin1996} or Theorem
3.4, page 289, in \cite{Conway1978})
to represent it as an infinite product over its zeros
\begin{eqnarray*}
Q(z)&=&\exp(c_1 z) \biggl(1+\frac{\i z}{\zeta_0^{+}}  \biggr)
\biggl(1+\frac{\i z}{\zeta_0^{-}}  \biggr)
\\
&&{}\times
\prod _{z_k \in{\mathfrak Z}}  \biggl(1+\frac{\i z}{z_k} \biggr)
\prod _{| n |\ge1}  \biggl(1+\frac{\i z}{\zeta
_n} \biggr)\exp
 \biggl(-\frac{\i z}{\zeta_n} \biggr)
\end{eqnarray*}
for some constant $c_1 \in{\mathbb C}$.
As the next step we rearrange the infinite product in the above formula
and obtain
%
%
\begin{eqnarray}\label{product_Q_zeta}
Q(z)&=&\exp(c_2 z) \biggl(1+\frac{\i z}{\zeta_0^{+}}  \biggr)
\biggl(1+\frac{\i z}{\zeta_0^{-}}  \biggr)
\nonumber\\[-8pt]\\[-8pt]
&&{}\times
\prod _{z_k \in{\mathfrak Z}}  \biggl(1+\frac{\i z}{z_k} \biggr)
\prod _{n \ge1}  \biggl(1+\frac{\i z}{\zeta_n} \biggr)
\biggl(1+\frac
{\i z}{\zeta_{-n}} \biggr)\nonumber
\end{eqnarray}
for some other constant $c_2 \in{\mathbb C}$, where the infinite product
converges absolutely since each term is $1+O(n^{-2})$ as $n\to\infty$.
Using definition of $P(z)$ (\ref{def_P}) and infinite product
representation of trigonometric functions (see formulas 1.431 in \cite
{Jeffrey2007}) we find that
%
%
\begin{equation}\label{product_P}
P(z)=\prod _{n \ge1}  \biggl(1+\frac{\i z}{n+\alpha}
\biggr)
\biggl(1+\frac{\i z}{-n+\alpha} \biggr).
\end{equation}
Combining equations (\ref{def_P}), (\ref{def_Q}),
(\ref{product_Q_zeta}) and (\ref{product_P}) we finally conclude that
for all $z \in{\mathbb C}$
%
%
\begin{equation}\label{proof_factn}
\frac{q}{q+\Psi(z)}=\frac{\exp(c_2 z)}{ (1+{\i z}/{\zeta
_0^{+}}  )  (1+{\i z}/{\zeta_0^{-}} )}
\prod _{z_k \in{\mathfrak Z}} \frac{1}{1+{\i z}/{z_k}}
\prod _{| n |\ge1} \frac{1+{\i z}/({n+\alpha
})}{1+
{\i z}/{\zeta_n}}.\hspace*{-28pt}
\end{equation}

First let us prove that $c_2=0$. Denote the left-hand side of
(\ref{proof_factn})
as $F_1(z)$ and right-hand side as $F_2(z)$.
Since $\Psi(z)$ is a characteristic exponent,
it must be $O(z^2)$ as $z \to\infty$, $z \in{\mathbb R}$, thus clearly
$z^{-1}\ln(F_1(z)) \to0$ as $z \to\infty$, $z\in{\mathbb R}$.
Using Lem\-ma~\ref{lemmma_asymp_for_product} we find that $z^{-1}\ln
(F_2(z)) \to c_2$ as $z \to\infty$ $z \in{\mathbb R}$, which implies that
$c_2=0$.

All that is left to do it to prove that ${\mathfrak Z}$ is an empty
set. The main tool is again Lemma \ref{lemmma_asymp_for_product}.
Assuming that $\sigma\ne0$ and using asymptotic expansion (\ref
{eq_asympt_roots1_1})
and Lem\-ma~\ref{lemmma_asymp_for_product} we find that the infinite
product in (\ref{proof_factn}) converges to a
constant as $z\to\infty$, $z\in{\mathbb R}$. Thus function $F_2(z)
\approx A_2
z^{-2-M}$ where $M$ is equal to the number of
elements in the set ${\mathfrak Z}$.
However, function $F_1(z) \approx A_1 z^{-2}$ as $z\to\infty$, $z\in
{\mathbb R}
$, thus $M=0$ and the set ${\mathfrak Z}$ must be empty.
In the case $\sigma\ne0$ the proof is identical,
except that both $F_1(z)$ and $F_2(z)$ behave like $Az^{-1}$, which can
be established by the asymptotic expression for
$\zeta_n$ given in
(\ref{eq_asympt_roots1_2})
and Lemma \ref{lemmma_asymp_for_product}.
\end{pf*}

Theorem \ref{thm_nu3_1} provides us with all the information about the
zeros of $q+\Psi(z)$ that we will need
later to prove results about Wiener--Hopf factors and perform numerical
computations. However, we can also compute explicitly
the sums of inverse powers of zeros. These results can be useful for
checking\vspace*{1pt} the accuracy, but more importantly, for approximating the
smallest solutions
$\zeta_0^{\pm}$. We assume that $\alpha\ne0$ and define for
$m\ge0$
%
%
\begin{equation}\label{def_Sm}
\Omega_m=\alpha^{-m-1}+ ( \zeta_0^{-}
)^{-m-1}+ (\zeta_0^{+} )^{-m-1}+\sum _{n\ge1}
 [ \zeta_n^{-m-1}+\zeta_{-n}^{-m-1}  ].
\end{equation}
Asymptotic expansions (\ref{eq_asympt_roots1_1}) and (\ref
{eq_asympt_roots1_2}) guarantee that the series converges absolutely
for $m\ge0$, thus the sequence $\{\Omega_m\}_{m\ge0}$ is correctly defined.
\begin{lemma}\label{Lemma_explicit_Sm}
The sequence $\{\Omega_m\}_{m\ge0}$ can be computed using the
following recurrence relation:
\[
\Omega_m=-\frac{1}{b_0} \Biggl[(m+1)b_{m+1}+\sum _{n=0}^{m-1}
\Omega
_n b_{m-n}  \Biggr], \qquad    m\ge0,
\]
where coefficients $\{b_n\}_{n\ge0}$ are defined as
%
%
\begin{eqnarray}\label{def_bn}
b_{2n}&=&\frac{(-1)^{n-1}\pi^{2n-1}}{(2n)!}  [ n(2n-1)\alpha
\sigma
^2+\pi^2\alpha(q+8n)-2n\gamma\rho ] \nonumber\\
b_{2n+1}&=&\frac{(-1)^{n}\pi^{2n}}{(2n+1)!}  \biggl[n(2n+1)\frac
{\gamma
\sigma^2}{\pi} -
\pi(4\pi^2 \alpha^2+4\gamma^2-\gamma q)\\
&&\hspace*{128.2pt}{}+\pi(2n+1)(4\gamma+\alpha
\rho)
\biggr].\nonumber
\end{eqnarray}
\end{lemma}
\begin{pf}
This statement is just an application of the following general result.
Assume that we have an entire function $H(z)$ which can be expressed
as an infinite product over the set of its zeros ${\mathfrak Z}$
\[
H(z)=\prod _{z_k \in{\mathfrak Z}} \biggl(1-\frac
{z}{z_k} \biggr).
\]
Taking derivative of $\ln(H(z))$ we find
\[
H'(z)=-H(z)\sum _{z_k \in{\mathfrak Z}} (z_k-z)^{-1}=-h(z)
\sum
 _{m\ge0}
\biggl[  \sum _{z_k \in{\mathfrak Z}} z_k^{-m-1}  \biggr] z^m,
\]
and the recurrence relation for $\sum_{z_k \in{\mathfrak Z}}
z_k^{-m-1}$ is obtained by expanding $H(z)$ and $H'(z)$ as a Maclaurin
series, multiplying two series in the right-hand side and comparing the
coefficients in front of $z^m$. The statement of Lemma \ref
{Lemma_explicit_Sm} follows by considering an entire function
\[
H(z)=q \pi(z-\alpha) Q(\i z),
\]
where $Q(z)$ is defined by (\ref{def_Q}). Function $H(z)$ has zeros at
$\{\alpha, \zeta_0^{\pm}, \zeta_n \}$, and one can check
that the Maclaurin expansion
is given by $H(z)=\sum_{n\ge0} b_n z^n$ where coefficients $b_n$ are
defined in (\ref{def_bn}).
\end{pf}

Finally we can state and prove our main results: expressions for
Wiener--Hopf factors and density of $S_{\tau}$.
\begin{theorem}\label{WH_factors_nu3} For $q>0$
%
%
\begin{eqnarray}\label{explicit_WH_factors}
\phi_q^{-}(z)&=&
\frac{1}{1+{\i z}/{\zeta_0^{+}}}
\prod _{n \ge1} \frac{1+{\i z}/({n+\alpha})}{1+{\i
z}/{\zeta_n}},\nonumber\\[-8pt]\\[-8pt]
\phi_q^{+}(z)&=&
\frac{1}{1+{\i z}/{\zeta_0^{-}}}
\prod _{n \le-1} \frac{1+{\i z}/({n+\alpha})}{1+{\i
z}/{\zeta_n}}.\nonumber
\end{eqnarray}
Infinite products converge uniformly on compact subsets of ${\mathbb
C}\setminus\i{\mathbb R}$.
The density of $S_{\tau}$ is given by
%
%
\begin{equation}\label{density_of_Mtauq}
\frac{\d}{\d x}   {\mathbb P}(S_{\tau} \le x)=- c_0^{-}\zeta
_0^{-
}e^{\zeta_0^{-}x}-
\sum _{k \le-1} c_k^{-} \zeta_{k} e^{\zeta_{k} x},
\end{equation}
where
%
%
\begin{eqnarray}\label{eq_ck_minus}
c_0^{-}&=&
\prod _{n \le-1} \frac{1-{\zeta_0^{-}}/({n+\alpha
})}{1-{\zeta_0^{-}}/{\zeta_n}},\nonumber\\[-8pt]\\[-8pt]
c_k^{-}&=&\frac{ 1-{\zeta_{k}}/({k+\alpha})}{1-{\zeta
_{k}}/{\zeta_0^{-}}}
\prod _{n \le-1,   n\ne k} \frac{1-{\zeta
_{k}}/({n+\alpha
})}{1-{\zeta_{k}}/{\zeta_n}}.\nonumber
\end{eqnarray}
\end{theorem}
\begin{pf}
Expressions (\ref{explicit_WH_factors}) for Wiener--Hopf factors are
obtained using factorization identity (\ref{eq_big_factorization})
and Lemmas \ref{Lemma1} and \ref{lemmma_asymp_for_product}. Expression
(\ref{density_of_Mtauq}) for the density of $S_{\tau}$ is
derived by computing the inverse Fourier transform via residues.
\end{pf}
\begin{remark}
Theorem \ref{WH_factors_nu3} remains true for $q=0$ if $\mu<0$. In this
case ${\mathbb E}X_1 <0$ and $S_{\tau}\to S_{\infty}$
and $I_{\tau} \to-\infty$ as $q\to0^+$. From the analytical point of
view we have $\zeta_0^{-}<0$ and
$\zeta_0^{+}=0$ [see Figure
\ref{fig_plot} and (\ref{eq_Psiq0})]. If $\mu=0$, then
${\mathbb E}
X_1=0$ and the process $X_t$ oscillates; thus $S_{\infty}=I_{\infty
}=\infty$, which is expressed analytically by the fact that $\zeta
_0^{+}=\zeta_0^{-}=0$.
\end{remark}

\section{A family of L\'evy processes}\label{section_results_nu4}

\begin{definition}
We define a $\beta$-family of L\'evy processes by the generating triple
$(\mu,\sigma,\nu)$, where the
density of the L\'evy measure is defined as
%
%
\begin{equation}\label{def_nu4}
\nu(x)=c_1\frac{e^{-\alpha_1 \beta_1 x}}{(1-e^{-\beta_1
x})^{\lambda
_1}}{\mathbf I}_{\{x>0\}}+c_2\frac{e^{\alpha_2 \beta_2
x}}{(1-e^{\beta
_2 x})^{\lambda_2}}{\mathbf I}_{\{x<0\}}
\end{equation}
and parameters satisfy $\alpha_i>0$, $\beta_i>0$, $c_i\ge0$ and
$\lambda_i \in(0,3)$. This L\'evy measure has exponential tails; thus
we will use
the cut-off function $h(x)\equiv x$ in (\ref{levy_khintchine}).
\end{definition}

The $\beta$-family is quite rich: in particular, by controlling
parameters $\lambda_i$,
we can obtain an arbitrary behavior of small jumps, and parameters
$\alpha_i$ and $\beta_i$ are responsible for the tails of the L\'evy measure
(which are always exponential). Parameters $c_i$ control the total
``intensity'' of positive/negative jumps. The processes in $\beta
$-family are
similar to the generalized tempered stable processes (see \cite{Cont})
which were also named KoBoL processes in \cite{BoLe} and \cite{Boyarchenko}
\[
\nu(x)=
c_{+} \frac{e^{-\alpha_{+} x}}{ x ^{\lambda_{+
}}}{\mathbf
I}_{\{x>0\}}+
c_{-} \frac{e^{\alpha_{-} x}}{| x |^{\lambda_{-
}}}{\mathbf I}_{\{x<0\}}.
\]
In fact we can obtain the above measure as the limit of L\'evy measures
in $\beta$-family. If we set $c_1=c_{+}\beta^{\lambda_{+}}$,
$c_2=c_{-}\beta^{\lambda_{-}}$, $\alpha_1=\alpha_{+
}\beta^{-1}$,
$\alpha_2=\alpha_{-}\beta^{-1}$, $\beta_1=\beta_2=\beta$
and let
$\beta\to0^+$ we see that the L\'evy measure defined in
(\ref{def_nu4}) will converge to the L\'evy measure of the generalized
tempered stable process.
Next, when $\lambda_1=\lambda_2$, the processes in $\beta$-family are
similar to the
tempered stable processes (see \cite{BaLe}). If we restrict the parameters
even further, $c_1=c_2$, $\lambda_1=\lambda_2$
and $\beta_1=\beta_2$ so that the small positive/negative jumps have
the same behavior, while large jumps which are
controlled by $\alpha_i$ may be different, and we obtain a process very
similar to the CGMY family defined in \cite{CGMY2002}.
Finally, if $c_i=4$, $\beta_i=1/2$, $\lambda_i=2$ and $\alpha
_1=1-\alpha
$ and $\alpha_2=1+\alpha$, we obtain the process $X_t$ discussed in
Section \ref{section_results_nu3}.

If we restrict parameters as $\sigma=0$, $\beta_1=\beta_2$, $\lambda
_1=\lambda_2$ (and $\mu$ uniquely specified in terms of these
parameters), then $\beta$-family reduces to
a family of Lamperti-stable processes, which can be obtained by
Lamperti transformation from the stable processes conditioned to stay positive
(see the original paper \cite{Lamperti1972} by Lamperti for the
definition of this transformation and its various properties).
Spectrally one-sided Lamperti-stable processes appeared in
\cite{BertoinYor2001} and \cite{Patie2009}, and two-sided processes were
studied in \cite{CaCha,Caballero2008,Chaumont2009} and \cite{KyPaRi}.
Lamperti-stable processes are a particularly interesting subclass of
the $\beta$-family since they offer many examples of fluctuation
identities related to Wiener--Hopf factorization which can be computed
in closed form
(see \cite{CaCha,Chaumont2009} and \cite{KyPaRi}).

In the following proposition we derive a formula for the characteristic
exponent $\Psi(z)$ for processes in the $\beta$-family.
As we will see, the characteristic exponent can be expressed in terms
of beta and digamma functions (see Chapter 6 in \cite{AbramowitzStegun}
or Section 8.3 in \cite{Jeffrey2007})
%
%
\begin{equation}\label{eq_def_beta_digamma}
\mathrm{B}(x;y)=\frac{\Gamma(x)\Gamma(y)}{\Gamma(x+y)},\qquad
\psi
(x)=\frac{\d}{\d x}\ln(\Gamma(x)),
\end{equation}
which justifies the name of the family.
\begin{proposition} If $\lambda_i \in(0,3)\setminus\{1,2\}$, then
%
%
\begin{eqnarray}\label{def_psi_nu4}
\Psi(z)&=&\frac{\sigma^2z^2}2+\i\rho z
-\frac{c_1}{\beta_1} \mathrm{B} \biggl(\alpha_1-\frac
{\i
z}{\beta_1}; 1-\lambda_1  \biggr)\nonumber\\[-8pt]\\[-8pt]
&&{}-
\frac{c_2}{\beta_2} \mathrm{B} \biggl(\alpha_2+\frac{\i z}{\beta_2};
1-\lambda
_2  \biggr)+\gamma,\nonumber
\end{eqnarray}
where
\begin{eqnarray*}
\gamma&=&\frac{c_1}{\beta_1} \mathrm{B} (\alpha_1; 1-\lambda_1
)+\frac{c_2}{\beta_2} \mathrm{B} (\alpha_2; 1-\lambda_2),\nonumber\\
\rho&=&\frac{c_1}{\beta_1^2} \mathrm{B} (\alpha_1; 1-\lambda_1
)\bigl(\psi(1+\alpha_1-\lambda_1)-\psi(\alpha_1)\bigr)
\\
&&{}-
\frac{c_2}{\beta_2^2} \mathrm{B} (\alpha_2; 1-\lambda_2
)\bigl(\psi
(1+\alpha_2-\lambda_2)-\psi(\alpha_2)\bigr)-\mu.
\end{eqnarray*}
If $\lambda_1$ or $\lambda_2\in\{1,2\}$ the characteristic exponent can
be computed using the following two integrals:
%
%
\begin{eqnarray}\label{int_nu_lambda1}\qquad
&& \int _0^{\infty} (e^{\i xy}-1-\i xy)\frac{e^{-\alpha\beta x}
}{1-e^{-\beta x}}\,\d x\nonumber\\[-8pt]\\[-8pt]
&&\qquad=
 -\frac{1}{\beta}  \biggl[ \psi \biggl(\alpha
-\frac{\i
y}{\beta} \biggr)-\psi(\alpha) \biggr]-\frac{\i y}{\beta^2} \psi
'(\alpha)
\nonumber\\
\label{int_nu_lambda2}
&& \int _0^{\infty} (e^{\i xy}-1-\i xy)\frac{e^{-\alpha\beta x}
}{(1-e^{-\beta x})^2}\,\d x
\nonumber\\[-8pt]\\[-8pt]
&&\qquad  =-\frac{1}{\beta} \biggl(1-\alpha+\frac{\i
y}{\beta
} \biggr) \biggl[ \psi \biggl(\alpha-\frac{\i y}{\beta} \biggr)-\psi(\alpha
)
\biggr]-\frac{\i y(1-\alpha)}{\beta^2} \psi'(\alpha).\nonumber
\end{eqnarray}
\end{proposition}
\begin{pf}
First we assume that $\lambda\in(0,1)$. Performing change of
variables $u=\exp(-\beta x)$ and using integral representation for beta
function (formula 8.380.1 in
\cite{Jeffrey2007}) we find that
\begin{eqnarray*}
&&\beta\int _{0}^{\infty}  (e^{\i xz}-1-\i xz )
\frac
{e^{-\alpha\beta x}}{(1-e^{-\beta x})^{\lambda}}\,\d x\\
&&\qquad=
\mathrm{B} \biggl( \alpha-\frac{\i z}{\beta};1-\lambda \biggr)-
\mathrm{B} ( \alpha;1-\lambda )-z  \biggl[\frac{\d}{\d z}
\mathrm{B} \biggl( \alpha-\frac{\i z}{\beta};1-\lambda
\biggr) \biggr]_{z=0},
\end{eqnarray*}
and we obtain the desired result (\ref{def_psi_nu4}). The left-hand
side of the above equation is analytic in $\lambda$ for
$\operatorname{Re}(\lambda)<3$, and the right-hand side is analytic
and well defined
for $\operatorname{Re}(\lambda)<3$, $\lambda\ne\{1,2\}$; thus by analytic
continuation they should be equal for $\lambda\in(0,3)\setminus\{
1,2\}$.

Assume that $\lambda=2$. Then using binomial series we can expand
\[
\bigl(1-\exp(-x)\bigr)^{-2}=\sum_{n\ge0} (n+1)\exp(-nx),
\]
which converges uniformly on $(\varepsilon,\infty)$
and obtain
\begin{eqnarray*}
&&\beta\int _{0}^{\infty}  (e^{\i xz}-1-\i xz )
\frac
{e^{-\alpha\beta x}}{(1-e^{-\beta x})^{2}} \,\d x
\nonumber\\[-1pt]
&&\qquad=
\sum _{n\ge0}  \biggl[\frac{n+1}{n+\alpha-{\i y}/{\beta
}}-\frac
{n+1}{n+\alpha}-\frac{\i y}{\beta}\frac{n+1}{(n+\alpha)^2}
\biggr]\\[-1pt]
&&\qquad= \biggl(1-\alpha+\frac{\i y}{\beta}  \biggr)\sum _{n\ge0}
\biggl[\frac{1}{n+\alpha-{\i y}/{\beta}}-\frac{1}{n+\alpha}  \biggr]\\[-1pt]
&&\qquad\quad{} - \frac{\i y}{\beta}(1-\alpha)\sum _{n\ge0} \frac
{1}{(n+\alpha)^2},
\end{eqnarray*}
and using the series representation for digamma function (formula
8.362.1 in \cite{Jeffrey2007}) we obtain (\ref{int_nu_lambda2}).
Derivation of formula (\ref{int_nu_lambda1})
corresponding to the case $\lambda=1$ is identical.
\end{pf}

The following theorem is the analogue of Theorem \ref{thm_nu3_1}, and
it is the main result in this section.
\begin{theorem}\label{thm_nu4_1} Assume that $q>0$ and that $\Psi(z)$
is given by (\ref{def_psi_nu4}).

\begin{longlist}
\item Equation $q+\Psi(\i\zeta)=0$ has infinitely many solutions,
all of which are real and simple.
They are located as follows:
%
%
\begin{eqnarray}
\zeta_0^{-} &\in&(-\beta_1\alpha_1,0), \nonumber\\[-1pt]
\zeta_0^{+} &\in& (0,\beta_2\alpha_2), \nonumber\\[-9pt]\\[-9pt]
\zeta_n &\in& \bigl(\beta_2(\alpha_2+n-1), \beta_2(\alpha_2+n)\bigr),\qquad
n\ge
1 ,\nonumber\\[-1pt]
\zeta_n &\in&\bigl(\beta_1(-\alpha_1+n),\beta_1(-\alpha_1+
n+1)\bigr),\qquad
n\le-1.\nonumber
\end{eqnarray}
\item If $\sigma\ne0$ we have
%
%
\begin{eqnarray}\label{asympt_zetan_nu41}
\zeta_{n+1}&=&\beta_2(n+\alpha_2)
+\frac{2 c_2}{\sigma^2 \beta_2^2\Gamma(\lambda_2)}
(n+\alpha_2)^{\lambda_2-3}\nonumber\\
&&{}+O(n^{\lambda_2-3-\varepsilon}),\qquad
n\to+\infty,\nonumber\\[-8pt]\\[-8pt]
\zeta_{n-1}&=&\beta_1(n-\alpha_1)
-\frac{2 c_1}{\sigma^2\beta_1^2\Gamma(\lambda_1)}
(-n+\alpha_1)^{\lambda_1-3}\nonumber\\
&&{}+
O( n^{\lambda_1-3-\varepsilon}),\qquad
n\to-\infty.\nonumber
\end{eqnarray}
\item
If $\sigma=0$ we have
%
%
\begin{eqnarray}\label{asympt_zetan_nu42}
\zeta_{n+\delta}&=&\beta_2(n+\alpha_2+\omega_0)
+A(n+\alpha_2+\omega_0)^{\lambda}\nonumber\\[-8pt]\\[-8pt]
&&{}+O(n^{\lambda
-\varepsilon
}),\qquad
n\to+\infty,\nonumber
\end{eqnarray}
where coefficients $w_0$ and $A$ are presented in Table \ref{tab1},
$\delta\in\{0,1\}$ depending on the
signs of $w_0$ and $A$ and
\[
x_0=\frac{1}{\pi} \arctan \biggl( \sin(\pi\lambda_2)  \biggl(\frac
{c_1\beta
_2^{\lambda_2} \Gamma(1-\lambda_1)}{c_2\beta_1^{\lambda_1}\Gamma
(1-\lambda_2)}- \cos(\pi\lambda_2) \biggr)^{-1}  \biggr).
\]
The corresponding results for $n \to-\infty$ can be obtained by
symmetry considerations.

%
%
\begin{table}[b]
\caption{Coefficients for asymptotic expansion of $\zeta_n$ when
$\sigma=0$}\label{tab1}
\begin{tabular*}{\tablewidth}{@{\extracolsep{\fill}}lccc@{}}
\hline
&$\bolds{\omega_0}$ & $\bolds{A}$ & $\bolds{\lambda}$ \\
\hline
$\lambda_1<2$, $\lambda_2<2$ & 0 &
$ \frac{c_2}{\rho\beta_2 \Gamma(\lambda_2)}$ &
$\lambda
_2-2$ \\[4pt]
$\lambda_1<2$, $\lambda_2>2$ & $2-\lambda_2$ &
$ -\frac{\sin(\pi\lambda_2) \beta_2^3 \rho}{\pi c_2
\Gamma(1-\lambda_2)}$ & $2-\lambda_2$ \\[4pt]
$\lambda_1>2$, $\lambda_2<\lambda_1$ &0 &
$ \frac{c_2\beta_1^{\lambda_1}}{c_1 \beta_2^{\lambda
_1-1}\Gamma(1-\lambda_1)\Gamma(\lambda_2)}$ &
$\lambda_2-\lambda_1$ \\[4pt]
$\lambda_1>2$, $\lambda_2>\lambda_1$ &
$2-\lambda_2$ &
$ -\frac{\sin(\pi\lambda_2)}{\pi} \frac{c_1\beta
_2^{\lambda_1+1}
\Gamma(1-\lambda_1)}{c_2\beta_1^{\lambda_1}\Gamma(1-\lambda_2)}$
& $\lambda_1-\lambda_2$ \\[4pt]
$\lambda_1>2$, $\lambda_2=\lambda_1$ & $x_0$ &
$ -\rho\frac{\sin(\pi x_0)^2}{\pi^2}\frac{\beta_2^3
}{c_2}\Gamma(\lambda_2)$ & $2-\lambda_2$ \\
\hline
\end{tabular*}
\end{table}

\item Function $q(q+\Psi(z))^{-1}$ can be factorized as follows:
%
%
\begin{eqnarray}\label{big_factorization2}\quad
\frac{q}{q+\Psi(z)}&=&\frac{1}{ (1+{\i z}/{\zeta_0^{+}}
 ) (1+{\i z}/{\zeta_0^{-}}  )}
\prod _{ n \ge1} \frac{1+{\i z}/({\beta_2(n-1+\alpha
_2)})}{1+{\i z}/{\zeta_n}}\nonumber\\[-8pt]\\[-8pt]
&&{}\times
\prod _{ n \le-1} \frac{1+{\i z}/({\beta_1(n+1-\alpha
_1)})}{1+{\i z}/{\zeta_n}},\nonumber
\end{eqnarray}
where the infinite products converge uniformly on the compact subsets
of the complex plane excluding zeros/poles of $q+\Psi(z)$.
\end{longlist}
\end{theorem}
\begin{remark}
When $\sigma=0$ the remaining cases $\lambda_1<2$, $\lambda_2=2$ and
$\lambda_1=2$, $0<\lambda_2<3$ are not covered
by Theorem \ref{thm_nu4_1}. The interested reader can derive these
asymptotic expansions by
using formulas (\ref{int_nu_lambda1}), (\ref{int_nu_lambda2}) and
the following results for the digamma function (see formulas 6.3.7 and
6.3.18 in \cite{AbramowitzStegun}):
\[
\psi(1-z)=\psi(z)+\pi\cot(\pi z),\qquad
\psi(z)=\ln(z)-\frac{1}{2z}+O(z^{-2}),   \qquad  z\to\infty.
\]
\end{remark}
\begin{pf*}{Proof of Theorem \ref{thm_nu4_1}}
The proof of (i) is very similar to the corresponding part of the proof
of Theorem \ref{thm_nu3_1}. We separate equation
$q+\Psi(\i\zeta)=0$ into a jump part and a diffusion part, find points
where the jump part goes to infinity and by
analyzing the signs we conclude that on every interval between these
points there should exist a solution.

The proof of (ii) and (iii) is based on the following two asymptotic
formulas as $\zeta\to+\infty$:
\begin{eqnarray*}
\mathrm{B}(\alpha+\zeta;\gamma)&=&\Gamma(\gamma)\zeta^{-\gamma
}
\biggl[1-\frac
{\gamma(2\alpha+\gamma-1)}{2\zeta}+O(\zeta^{-2}) \biggr],\\
\mathrm{B}(\alpha-\zeta;\gamma)&=&\Gamma(\gamma)\frac{\sin(\pi
(\zeta-\alpha
-\gamma))}{\sin(\pi(\zeta-\alpha))}
\zeta^{-\gamma}  \biggl[1+\frac{\gamma(2\alpha+\gamma-1)}{2\zeta
}+O(\zeta
^{-2}) \biggr].
\end{eqnarray*}
The first asymptotic expansion follows from the definition of beta
function (\ref{eq_def_beta_digamma}) and formula 6.1.47 in
\cite{AbramowitzStegun}, while the second formula
can be reduced to the first one by applying a reflection formula for
the gamma function.

If $\sigma\ne0$ and $\zeta\to+\infty$ we use (\ref{def_psi_nu4}) and
the above formulas and rewrite equation $q+\Psi(\i\zeta)=0$ as
\begin{eqnarray*}
\frac{\sin ( \pi ({\zeta}/{\beta_2}-\alpha
_2+\lambda_2
) )}
{\sin ( \pi ({\zeta}/{\beta_2}-\alpha_2
) )}&=&
\frac{\sigma^2\beta_2^{\lambda_2}}{2c_2\Gamma(1-\lambda_2)} \zeta
^{3-\lambda_2}+O(\zeta^{1-\lambda_2})\\
&&{}+O(\zeta^{\lambda_1-\lambda_2}),
\end{eqnarray*}
while if $\sigma=0$ and $\zeta\to+\infty$ we have
\begin{eqnarray*}
\frac{\sin ( \pi ({\zeta}/{\beta_2}-\alpha
_2+\lambda
_2 ) )}
{\sin ( \pi ({\zeta}/{\beta_2}-\alpha_2
) )} &=&
\frac{\beta_2^{\lambda_2}\rho}{c_2\Gamma(1-\lambda_2)}\zeta
^{2-\lambda_2} +
\frac{c_1\beta_2^{\lambda_2} \Gamma(1-\lambda_1)}{c_2\beta
_1^{\lambda
_1}\Gamma(1-\lambda_2)}\zeta^{\lambda_1-\lambda_2}\\
&&{} + O(\zeta
^{1-\lambda
_2})+O(\zeta^{\lambda_1-\lambda_2-1}).
\end{eqnarray*}
Asymptotic expansions (\ref{asympt_zetan_nu41}) and (\ref
{asympt_zetan_nu42}) can be derived from the above formulas using
the same method as in the proof of Theorem \ref{thm_nu3_1}.\vadjust{\goodbreak}

In order to prove factorization identity (\ref{big_factorization2}) and
the fact that there are no other roots, we use exactly the same approach
as in the proof of Theorem \ref{thm_nu3_1}. Again we choose an entire
function $P(z)$ which has
zeros at the poles of $\Psi(z)$ with the same multiplicity, and the
choice is obvious due to (\ref{def_psi_nu4}):
%
%
\begin{equation}\label{def_PQ2}
P(z)=
\biggl[\Gamma \biggl(\alpha_1-\frac{\i z}{\beta_1} \biggr)\Gamma
\biggl(\alpha
_2+\frac{\i z}{\beta_2} \biggr) \biggr]^{-1},
\end{equation}
and function $Q(z)$ is defined by as $q^{-1}(1+\Psi(z))P(z)$.
Function $P(z)$ can be expanded in infinite product using Euler's
formula (see formula 6.1.3 in \cite{AbramowitzStegun}).
Using (\ref{def_psi_nu4}) and (\ref{def_PQ2}) and asymptotics for gamma
function we find that function $Q(z)$ has order equal
to one, and thus again we can use Hadamard's theorem to expand it as
infinite product,
and finally we use asymptotics for infinite products supplied by Lemma
\ref{lemmma_asymp_for_product} and asymptotics for $\zeta_n$ given by
(\ref{asympt_zetan_nu41}) and (\ref{asympt_zetan_nu41}) to prove factorization
identity (\ref{big_factorization2}).
\end{pf*}

We can also derive a result similar to Lemma \ref{Lemma_explicit_Sm}
using the entire function $Q(z)$ defined by
(\ref{def_PQ2}). While there is no closed form expression for
derivatives
of gamma function, they can be easily computed numerically. Our final
result in this section is the analogue of
Theorem \ref{WH_factors_nu3}, and the proof is identical.
\begin{theorem}\label{WH_factors_nu4} For $q>0$
\begin{eqnarray*}
\phi_q^{-}(z)&=&
\frac{1}{1+{\i z}/{\zeta_0^{+}}}
\prod _{n \ge1} \frac{1+{\i z}/({\beta_2(n-1+\alpha
_2)})}{1+{\i z}/{\zeta_n}},\\
\phi_q^{+}(z)&=&
\frac{1}{1+{\i z}/{\zeta_0^{-}}}
\prod _{n \le-1} \frac{1+{\i z}/({\beta_1(n+1-\alpha
_1)})}{1+{\i z}/{\zeta_n}}.
\end{eqnarray*}
Infinite products converge uniformly on compact subsets of ${\mathbb
C}\setminus\i{\mathbb R}$.
The density of $S_{\tau}$ is given by
\[
\frac{\d}{\d x}   {\mathbb P}(S_{\tau}\le x)=- c_0^{-}\zeta
_0^{-
}e^{\zeta_0^{-}x}-
\sum _{k \le-1} c_k^{-} \zeta_{k} e^{\zeta_{k} x},
\]
where
\begin{eqnarray*}
c_0^{-}&=&
\prod _{n \le-1} \frac{1-{\zeta_0^{-}}/({\beta
_1(n+1-\alpha_1)})}{1-{\zeta_0^{-}}/{\zeta_n}},\\
c_k^{-}&=&\frac{ 1-{\zeta_{k}}/({\beta_1(k+1-\alpha
_1)})}{1-
{\zeta_{k}}/{\zeta_0^{-}}}
\prod _{n \le-1,   n\ne k} \frac{1-{\zeta_{k}}/({\beta
_1(n+1-\alpha_1)})}{1-{\zeta_{k}}/{\zeta_n}}.
\end{eqnarray*}
\end{theorem}
\begin{remark}
In the case $\sigma=0$, $\lambda_i<2$, and $\rho>0$ we have a process of
bounded variation and negative drift, and thus the distribution of
$S_{\tau}$ will have an atom at zero, which can be computed using the
following formula:
\begin{eqnarray*}
{\mathbb P}(S_{\tau}=0)&=&\lim _{z\to+\infty} {\mathbb E}
[ e^{-z
S_{\tau}}
 ]=\lim _{z\to+\infty} \phi_q^{+}(\i z)\\
 &=&\frac
{-\zeta
_0^{-}}{\alpha_1 \beta_1} \prod _{n\le-1} \frac{\zeta
_n}{\beta_1(n-\alpha_1)}.
\end{eqnarray*}
Using asymptotic relation (\ref{asympt_zetan_nu42}) one can see that
the above infinite product converges to
a number between zero and one.
\end{remark}

\section{Implementation and numerical results}\label{section_implementation}

In this section we discuss implementation details for computing the
probability density function of $S_{\tau(q)}$ and $S_t$.
In order to illustrate the main ideas we will use the process $X_t$ defined
in Section~\ref{section_results_nu3}; however, the implementation for a
general $X_t$ from the $\beta$-family would be quite
similar.
Our main tools are Theorem \ref{WH_factors_nu3} and asymptotic
expansion for $\zeta_n$ given in Theorem \ref{thm_nu3_1}.

First let us discuss the computation of density of $S_{\tau}$. The
first step would be to compute solutions to
equation $q+\Psi(\i\zeta)=0$, and for $q$ real this is a simple task:
for $n$ large we use Newton's method which is started from the
point
given by asymptotic expansion (\ref{eq_asympt_roots1_1}) or (\ref
{eq_asympt_roots1_2}). To compute $\zeta_0^{\pm}$ or $\zeta
_n$ with
$n$ small we use localization result (\ref{eq_localization}) and the
secant (or bisection) method to get the starting point for Newton's
iteration. Overall this part of the algorithm is very computationally
efficient and can be made even faster if
we compute different $\zeta_n$ in parallel.

The second step is to compute coefficient $c_k^{-}$ which are
given by (\ref{eq_ck_minus}). Each term in the infinite product is
$1+O(n^{-2})$; however, as we show in Proposition \ref{prop_conv_accel},
we can considerably improve
convergence by using our knowledge of the asymptotic expansion for
$\zeta_n$.
The final step is to compute the density of $S_{\tau}$ using
formula~(\ref{density_of_Mtauq}). Note that the series converges exponentially
for $x>0$. When $x$ is small the convergence is slow, and the
asymptotic behavior as $x\to0^+$
would depend on the decay rate of coefficient $c_k^{-}$; however,
we were unable to prove any results in this direction.
\begin{proposition}\label{prop_conv_accel}
Assume that $\zeta_n=n+\beta+\frac{A_1}{(n+\beta)}+\frac
{A_2}{(n+\beta
)^2}+O(n^{-3})$ as $n\to+\infty$.
Then as $N\to+\infty$ we have
%
%
\begin{eqnarray}\label{eq_conv_imprv}
&&
\prod _{n\ge N} \frac{1+{z}/({n+\alpha})}{1+
{z}/{\zeta_n}}\nonumber\\
&&\qquad=
\frac{\Gamma(N+\alpha)\Gamma(N+\beta+z)}{\Gamma(N+\beta)\Gamma
(N+\alpha+z)}
\nonumber\\[-8pt]\\[-8pt]
&&\qquad\quad{}\times
\exp \bigl[A_1 \bigl(f_{1,1}(\beta,\beta;N)-f_{1,1}(z+\beta,\beta
;N) \bigr)
\nonumber\\
&&\qquad\quad\hspace*{30.3pt}{} + A_2 \bigl(f_{1,2}(\beta,\beta;N)-f_{1,2}(z+\beta
,\beta ;N) \bigr)+O(N^{-3})
 \bigr],\nonumber
\end{eqnarray}
where $f_{\alpha_1,\alpha_2}(z_1,z_2;N)$ can be computed as
follows:
%
%
\begin{eqnarray}\label{asympt_faazzN}
&&
f_{\alpha_1,\alpha_2}(z_1,z_2;N)\nonumber\\
&&\qquad=\sum _{k\ge0}
\frac{{-\alpha_2\choose k}(z_2-z_1)^k}{\alpha_1+\alpha
_2+k-1}(z_1+N)^{1-\alpha_1-\alpha_2-k}\nonumber\\[-8pt]\\[-8pt]
&&\qquad\quad{}+(z_1+N)^{-\alpha_1}
(z_2+N)^{-\alpha_2} \biggl[\frac12+\frac{\alpha_1}{12(z_1+N)}+\frac
{\alpha_2}{12(z_2+N)}\biggr]
\nonumber\\
&&\qquad\quad{}+O(N^{-\alpha_1-\alpha_2-3}).\nonumber
\end{eqnarray}
\end{proposition}
\begin{pf}
First we define for $\alpha_1+\alpha_2>1$
\[
f_{\alpha_1,\alpha_2}(z_1,z_2;N)=\sum _{ n \ge N}
(n+z_1)^{-\alpha
_1} (n+z_2)^{-\alpha_2}.
\]
The proof of the asymptotic expansion (\ref{asympt_faazzN}) is based on
the Euler--Maclaurin formula
\[
\sum _{n\ge N} f(n)=\int _N^{\infty} f(x)\,\d x + \frac
{f(N)}2-\frac{f'(N)}{12}+O\bigl(f^{(3)}(N)\bigr),
\]
where we take $f(x)=(x+z_1)^{-\alpha_1} (x+z_2)^{-\alpha_2}$. To obtain
(\ref{asympt_faazzN}) we compute the integral by changing variables
$y=(x+z_1)^{-1}$, expanding the resulting integrand in Taylor's series
at $y=0$ and integrating term by term.

In order to obtain formula (\ref{eq_conv_imprv}) we follow the steps of
the proof of Lemma \ref{lemmma_asymp_for_product}
\begin{eqnarray*}
\prod _{n\ge N} \frac{1+{z}/({n+\alpha})}{1+
{z}/{\zeta_n}}&=&
\prod _{n\ge N} \frac{1+{z}/({n+\alpha})}{1+
{z}/({n+\beta})}
\prod _{n\ge N} \frac{1+{z}/({n+\beta})}
{1+{z}/{\zeta_n}}\\
&=&\frac{\Gamma(N+\alpha)\Gamma(N+\beta+z)}{\Gamma(N+\beta
)\Gamma
(N+\alpha+z)} \prod _{n\ge N} \frac{1+{z}/({n+\beta})}
{1+{z}/{\zeta_n}}.
\end{eqnarray*}

Next, we approximate $\ln(1+\omega)=\omega+O(\omega^2)$ as $\omega
\to
0$ and obtain
\begin{eqnarray*}
&&
\sum _{n\ge N} \ln \biggl(\frac{1+{z}/({n+\beta})}
{1+{z}/{\zeta_n}} \biggr)\\
&&\qquad=\sum _{n\ge N} \ln \biggl(\frac
{\zeta
_n}{n+\beta} \biggr)-
\sum _{n\ge N} \ln \biggl(\frac{z+\zeta_n}
{z+n+\beta} \biggr)\\
&&\qquad= \sum _{n\ge N}  \biggl(\frac{A_1}{(n+\beta)^2}+\frac
{A_2}{(n+\beta)^3}+O(n^{-4}) \biggr) \\
&&\qquad\quad{}-
\sum _{n\ge N}  \biggl(\frac{A_1}{(z+n+\beta)(n+\beta)}+\frac
{A_2}{(z+n+\beta)(n+\beta)^2}+O(n^{-4}) \biggr),
\end{eqnarray*}
which completes the proof.
\end{pf}

Computing the density $p_t(x)=\frac{\d}{\d x}{\mathbb P}(S_t\le x)$
of $S_t$ at
a deterministic time $t > 0$ requires more work.
Our starting point is the the fact that the density
of $S_{\tau}$ [which we denote by $p^S(q,x)$] is the Laplace transform
of $q \times p_t(x)$
\begin{eqnarray*}
p^S(q,x)&=&\frac{\d}{\d x} {\mathbb P}\bigl(S_{\tau(q)} \le x\bigr)= \frac{\d
}{\d x}
\int _0^{\infty} q e^{-qt} {\mathbb P}(S_t\le x)\,\d t\\
&=& q\int _0^{\infty} e^{-qt} p_t(x) \,\d t.
\end{eqnarray*}
Thus $p_t(x)$ can be recovered as the following cosine transform:
\[
p_t(x)=
\frac{2}{\pi}e^{q_0t}\int _0^{\infty} \operatorname
{Re} \biggl[ \frac
{p^S(q_0+\i
u,x)}{q_0+\i u}  \biggr] \cos(tu)\,\d u, \qquad    q_0>0.
\]
We see that to compute this Fourier integral numerically we
need to be able to compute $p^S(q,x)$ for $q$ lying in some interval
$q\in[q_0,q_0+\i u_0]$ in the complex plane. The main problem is that
we need to solve many equations (\ref{eq_Psiq0}) with complex $q$.
While the asymptotic expansions for $\zeta_n$ presented in Theorem
\ref{thm_nu3_1}
are still true, we do not have any localization results in the complex plane.
It is certainly possible to compute the roots using the argument
principle, and originally all the computations were done by the author
using this method. However as we will see, there is a much more
efficient algorithm.

We need to compute the solutions of equation $q+\Psi(\i\zeta)=0$ for
all $q\in[q_0,q_0+\i u_0]$.
First we compute the initial values: the roots $\zeta_0^{\pm}$,
$\zeta_n$ for real value of $q=q_0$ using the method discussed above.
Next we consider each root as an implicit function of $u\dvtx\zeta_n(u)$
is defined as
\[
\Psi(\i\zeta_n(u))+(q_0+\i u)=0, \qquad \zeta_n(0)=\zeta_n.
\]
Using implicit differentiation we obtain a first order differential equation
\[
\frac{\d\zeta_n(u)}{\d u}=-\frac{1}{\Psi'(\i\zeta_n(u))}
\]
with initial condition $\zeta_n(0)=\zeta_n$. We compute the solution to
this ODE using a numerical scheme, for example, an adaptive Runge--Kutta
method, and at each step we correct the solution by applying several
iterations of Newton's method. Again, for different $n$ we can compute
$\zeta_n(u)$ in parallel.

%
%
\begin{figure}

\includegraphics{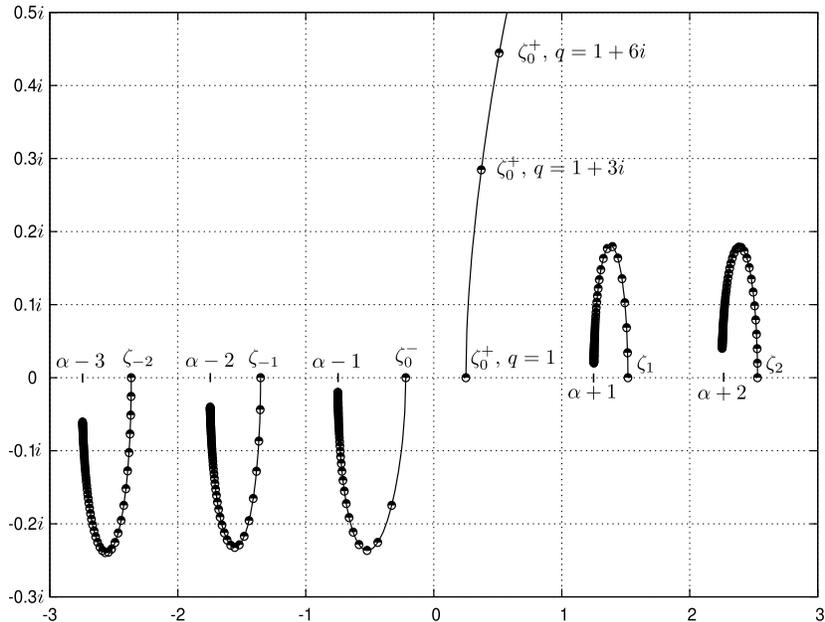}

\caption{The values of $\zeta_0^{\pm}$ and $\zeta_n$
for $q\in[1,1+200\i]$.}\label{fig_plot_cmplx}
\end{figure}

Figure \ref{fig_plot_cmplx} shows the result of this procedure. We have
used the following values of parameters: $\sigma=1$, $\mu=-0.1$ and
$\alpha=0.25$
and computed zeros $\zeta_0^{\pm}$ and $\zeta_n$ for $q\in
[1,1+200\i]$. The graph shows interesting
qualitative behavior: all zeros except $\zeta_0^{+}$ converge to
the closest pole of
$\Psi(\i\zeta)$ at $\alpha+n$, while $\zeta_0^{+}$ has no
pole nearby
[since $\Psi(\i\zeta)$ is regular at $\zeta=\alpha$] and it converges
to $\infty$ while always
staying in ${\mathbb C}^{+}$. If $\alpha<0$ the situation is exactly the
same, except that now $\zeta_0^{-}$ escapes to $\infty$
while always staying
in ${\mathbb C}^{-}$. We have repeated this procedure for many different
values of parameters, and from this numerical evidence
we can make some observations/conjectures.
It appears to be true that the roots never collide which means that we
have no
higher order solutions to $q+\Psi(z)=0$ for all $q\in{\mathbb C}$.
It also seems that the roots never cross the real line. All these
observations are based on numerical evidence, and we did not pursue
this any further to obtain rigorous proofs. However there is one fact
that we can prove rigorously: there are not going to appear any new,
unaccounted zeros. This could be proved
by an argument that we have used in the proof of Theorem
\ref{thm_nu3_1} to show that there are no extra zeros.

\begin{figure}

\includegraphics{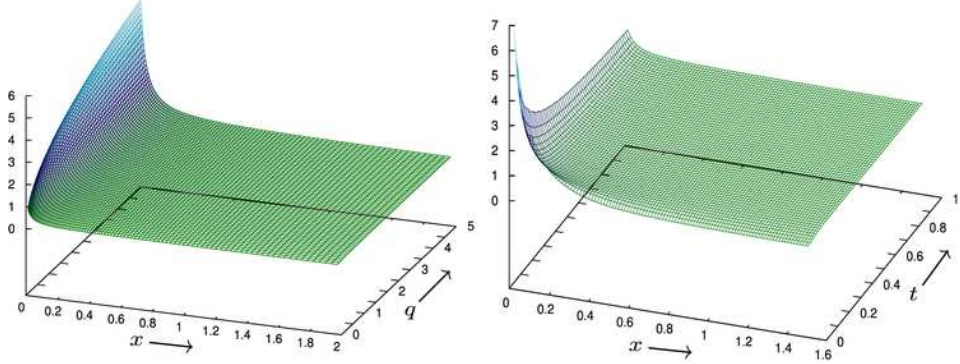}

\caption{Surface plot of $p^S(q,x)$ (left) and $p_t(x)$
(right).}\label{fig_surf}
\end{figure}

The results of our computations are presented in Figure \ref{fig_surf}.
The parameters are $\sigma=1$, $\mu=-0.1$ and $\alpha=0.25$;
the surface $p^S(q,x)$ is on the left and $p_t(x)$ on the right.

\section{Conclusion}\label{section_conclusion}

In this paper we have introduced a ten-parameter family of L\'evy
processes characterized by the fact that the
characteristic exponent is
a meromorphic function expressed in terms of beta and digamma functions.
This family is quite rich, and, in particular, it includes processes
with the complete range of
behavior of small positive/negative jumps. We have presented results of
the Wiener--Hopf factorization
for these processes, including semi-explicit formulas for Wiener--Hopf
factors and the density of the supremum process $S_{\tau}$.

%
%

These L\'evy processes might be used for modeling purposes whenever one
needs to compute distributions related to such functionals
as the first passage time, overshoot, extrema, last time before
achieving extrema, etc.
Some possible applications in Mathematical Finance and Insurance
Mathematics include pricing barrier, lookback and perpetual American options,
building structural models with jumps in Credit Risk, computing ruin
probabilities, etc.

Finally, we would like to mention that one can use the methods
presented in this paper, and, in particular, the technique to
solve $q+\Psi(z)=0$ for complex values of $q$ discussed in Section
\ref{section_implementation}, to compute
Wiener--Hopf factors arising from a more general factorization identity
(see Theorem 6.16 in \cite{Kyprianou})
\[
\frac{q}{q-\i w+\Psi(z)}=\Psi_q^{+}(w,z)\Psi_q^{-}(w,z),
\qquad  w,z \in{\mathbb R},
\]
where
\[
\Psi_q^{+}(w,z)=E  [e^{\i w \overline G_{\tau}+\i zS_{\tau
}} ],\qquad
\Psi_q^{-}(w,z)=E  [e^{\i w \underline{G}_{\tau}+\i
zI_{\tau
}} ],
\]
and $\overline G_t$ ($\underline G_t$) are defined as the last time
before $t$ when maximum (minimum) was achieved
\[
\overline G_t=\sup\{0\le u \le t \dvtx X_u=S_u \},  \qquad   \underline
G_t=\sup\{0\le u \le t \dvtx X_u=I_u \}.
\]

\section*{Acknowledgment}
The author would like to thank Andreas Kyprianou and Geert Van Damme
for many detailed comments and constructive suggestions.

%

%
\printaddresses

\end{document}